\documentclass{amsart}
\usepackage{amssymb}

\usepackage{amscd}

\newtheorem{theorem}{Theorem}[section]
\newtheorem{corollary}[theorem]{Corollary}

\newtheorem{lemma}[theorem]{Lemma}
\newtheorem{proposition}[theorem]{Proposition}

\theoremstyle{definition}

\theoremstyle{remark}

\begin{document}
\title{Deformations of the trivial line bundle and vanishing theorems}
\author{Herbert Clemens}
\author{Christopher Hacon}
\address{Mathematics Department, University of Utah}
\email{clemens@math.utah.edu\\
hacon@newmath.ucr.edu}
\date{November, 2000}
\maketitle

\begin{abstract}
This paper reproves a general form of the Green-Lazarsfeld ``generic
vanishing'' theorem and more recent strengthenings, as well as giving some
new applications.
\end{abstract}

\section{Introduction\protect\footnote{%
First author partially supported by NSF grant DMS-9970412}}

Our purpose is to generalize and give several applications of Theorem 3.2 of 
\cite{GL2} which characterizes the stalk at zero of the higher direct-image
sheaves 
\begin{equation}
R^{q}\pi _{*}\left( L\right)  \label{0.0}
\end{equation}
where $X_{0}$ is a compact K\"{a}hler manifold, $\Delta $ is the universal
cover of $\mathrm{Pic}^{0}\left( X_{0}\right) $, $\pi :X_{0}\times \Delta
\rightarrow \Delta $ is the standard projection, and 
\begin{equation*}
L\rightarrow X_{0}\times \Delta
\end{equation*}
is the (pullback of the) Poincare bundle. In \S 6 of \cite{GL2} Green and
Lazarsfeld also give a generalization to the case of 
\begin{equation}
R^{q}\pi _{*}\left( \Omega _{X_{0}}^{p}\otimes L\right) .  \label{0.1}
\end{equation}
A more complete picture of the Green-Lazarsfeld result emerged with the work
of Carlos Simpson \cite{S2} which explointed the associated Higgs bundle
structures to explain linearity and rationality of the strata defined by the
semicontinuous function given by the ranks of the geometric fibers of the
sheaves $\left( \ref{0.1}\right) .$ This very beautiful story has not been
brought together in a single, unified treatment, and such a treatment is
material to the generalization presented here. Therefore we ask the reader's
indulgence while we begin with a self-contained treatment of Theorem 6.1 of 
\cite{GL2}, a treatment which emphasizes the role of harmonic theory and
thereby helps clarify how Simpson's approach enters. We of course borrow
heavily from Green-Lazarsfeld \cite{GL1} and \cite{GL2} and Carlos Simpson's
work \cite{S1} and \cite{S2}.

It is then natural to generalize to the case of ``twisted coefficients.''
That is, in the spirit of the proof of Kawamata-Viehweg vanishing, we take a
cyclic covering 
\begin{equation*}
Y_{0}\rightarrow X_{0}
\end{equation*}
with simple-normal-crossing branch locus and push the above machinery on $%
Y_{0}$ forward to $X_{0}$ while keeping track of eigenspaces. This same
technique has been utilized by C. Mourougane to produce related results (see 
\cite{Mo}$)$. See also \cite{EV1} and \cite{Du} for related results and
techniques.

We will give several applications of this generalization. In \S 9, we give
an application of these methods to the zero sets of sections of a line
bundle on an abelian variety. In particular we prove a result on the
divisibility of linear series on abelian varieties analogous to a well known
theorem of Esnault and Viehweg concerning the zero sets of polynomials
(Theorem 2 of \cite{EV2}). Our result is sharp and in particular it implies
Theorem 1 of \cite{H} on the singularities of divisors on a principally
polarized abelian variety $(A,\Theta )$.

To indicate the utility these vanishing results, we indicate one simple
application. A result of Ein and Lazarsfeld (\cite{EL}, Proposition 3.5)
states that, for any effective divisor $D\in |N\Theta |$, the pair $%
(A,(1/N)D)$ is log canonical. In particular, for all points $p\in A,$%
\begin{equation*}
mult_{p}D\leq N\cdot \dim A.
\end{equation*}
Moreover, they show in Theorem 3.3 of \cite{EL} that, if $\Theta $ is
irreducible, then $(A,\Theta )$ is log terminal. As a corollary they show
that, if for $p\in A$ 
\begin{equation*}
mult_{p}(\Theta )\geq \dim A,
\end{equation*}
then equality holds and $\Theta $ must be reducible, so that there exists a
decomposition of principally polarized abelian varieties 
\begin{equation*}
(A,\Theta )\cong (A_{0},\Theta _{0})\times (A_{1},\Theta _{1}).
\end{equation*}
Since 
\begin{equation*}
mult_{p}(\Theta )=mult_{p_{0}(p)}\Theta _{0}+mult_{p_{1}(p)}\Theta _{1}
\end{equation*}
they repeat the argument and give another proof of a theorem of Smith and
Varley \cite{SV} that $A$ splits as a product of elliptic curves.

Corollary \ref{c9.4} of this paper complements the
Ein-Lazarsfeld-Smith-Varley results as follows. Suppose that $D\in \left|
N\Theta \right| $. Write 
\begin{equation*}
D=\sum\nolimits_{j}\left( k_{j}N+b_{j}^{\prime }\right) D_{j}
\end{equation*}
with $b_{j}^{\prime }<N$. Then 
\begin{equation}
D-N\left( \sum\nolimits_{j}k_{j}T_{j}\right) =\sum\nolimits_{j}b_{j}^{\prime
}D_{j}  \label{intdiv}
\end{equation}
is effective. Since 
\begin{equation*}
D-N\left( \sum\nolimits_{j}k_{j}T_{j}\right) \equiv N\left( \Theta
-\sum\nolimits_{j}k_{j}T_{j}\right)
\end{equation*}
the integral divisor 
\begin{equation}
\Theta -\sum\nolimits_{j}k_{j}T_{j}  \label{intdiv2}
\end{equation}
is semipositive. Writing $\Theta =(\Theta
-\sum\nolimits_{j}k_{j}T_{j})+\sum\nolimits_{j}k_{j}T_{j}$, a sum of
effective divisors, one sees by the Decomposition Theorem that all $k_{j}$
are either $0$ or $1$ and there is an isomorphism of principally polarized
abelian varieties 
\begin{equation*}
\left( A,\Theta \right) =\left( A^{\prime },\Theta ^{\prime }\right) \times
\left( A^{\prime \prime },\Theta ^{\prime \prime }\right)
\end{equation*}
with $\Theta -\sum\nolimits_{j}k_{j}T_{j}$ and $\sum\nolimits_{j}k_{j}T_{j}$
given by the pullbacks of $\Theta ^{\prime }$ and $\Theta ^{\prime \prime }$
respectively. Let $D^{\prime }$ be a divisor on $A^{\prime }$ such that $%
\sum b_{j}^{\prime }D_{j}$ is the pull back of $D^{\prime }$. And we
conclude:

\begin{theorem}
i) Suppose the principally polarized abelian variety $\left( A,\Theta
\right) $ is not a product of two principally polarized abelian varieties.
Then, 
\begin{equation*}
\left( A,\frac{1}{N}D\right) 
\end{equation*}
is log terminal.

ii) Suppose the principally polarized abelian variety $(A,\Theta )$ is not
the product of (principally polarized) elliptic curves. Then for all points $%
p\in A$ and $D\in |N\Theta |$, 
\begin{equation*}
mult_{p}D<N\cdot \dim A.
\end{equation*}
\end{theorem}

\begin{proof}
i) We have seen that, in the above discussion that 
\begin{equation*}
(A,\frac{1}{N}D)=(A^{\prime },\frac{1}{N}D^{\prime })\times (A^{\prime
\prime },\Theta ^{\prime \prime }).
\end{equation*}
If $A^{\prime \prime }=0$, $\left( A^{\prime },\frac{1}{N}D^{\prime }\right) 
$ Corollary \ref{c9.4} gives log terminality. If $A^{\prime }=0,$ $\Theta
^{\prime \prime }$ is irreducible and so by the result of Ein and Lazarsfeld 
\begin{equation*}
\left( A^{\prime \prime },\Theta ^{\prime \prime }\right)
\end{equation*}
is log terminal. So in either case $(A,\frac{1}{N}D)$ is log terminal.

ii) Again we write 
\begin{equation*}
(A,\frac{1}{N}D)=(A^{\prime },\frac{1}{N}D^{\prime })\times (A^{\prime
\prime },\Theta ^{\prime \prime })
\end{equation*}
with $(A^{\prime },\frac{1}{N}D^{\prime })$ log terminal. So in particular,
for all points $p^{\prime }\in A^{\prime }$, 
\begin{equation*}
mult_{p^{\prime }}D^{\prime }<N\cdot \dim A^{\prime }.
\end{equation*}
Suppose now that $mult_{p}D=N\cdot \dim A$. Then $A^{\prime }=0$, $D=N\Theta
^{\prime \prime }=N\Theta $ and $mult_{p}\Theta =\dim A$. Therefore by the
result of Smith and Varley, $(A,\Theta )$ must split as a product of
(principally polarized) elliptic curves.
\end{proof}

In \S 10, we turn our attention to the properties of the pushforwards of
dualizing sheaves. For a surjective morphism of projective varieties 
\begin{equation*}
f:X\longrightarrow X^{\prime },
\end{equation*}
$X$ smooth, we extend the Green-Lazarsfeld generic vanishing Theorem \ref
{c6.1} to (twists of) pushforwards of the dualizing sheaf on $X$.

Ein and Lazarsfeld \cite{EL} have recently shown that this type of generic
vanishing theorem can be used very effectively in the study of the geometry
of varieties of maximal Albanese dimension (that is, varieties whose
Albanese map is generically immersive). For irregular varieties not of
maximal Albanese dimension, the results of \S 10 are needed to study the
pushforwards $\omega _{X}$ under the Albanese map. For example, in \cite{CH}%
, these results are applied to the proof of an effective version of a
conjecture of Ueno. Namely, if $\kappa (X)=0$ and $h^{0}(\omega _{X})=1$
then ${alb_{X}}_{*}(\omega _{X})=\mathcal{O}_{Alb(X)}$. In particular for a
generic fiber $F$ one has $h^{0}(\omega _{F})=1$.

\section{Deformations of line bundles\label{2}}

Let 
\begin{equation}
L\overset{p}{\longrightarrow }X\overset{\pi }{\longrightarrow }\Delta
\label{2.2}
\end{equation}
be a holomorphic line bundle over the total space of a deformation $X/\Delta 
$ of complex manifolds. We view sections of $L$ as functions 
\begin{equation*}
f:L^{\vee }\rightarrow \Bbb{C}
\end{equation*}
such that 
\begin{equation*}
\left[ \chi ,f\right] =f
\end{equation*}
where $\chi $ is the holomorphic Euler vector field on $L^{\vee }$
associated to the $\Bbb{C}^{*}$-action.

We first claim that, given a transversely holomorphic trivialization 
\begin{equation*}
F_{\sigma }:X\overset{\left( \sigma ,\pi \right) }{\longrightarrow }
X_{0}\times \Delta ,
\end{equation*}
we can make compatible trivializations 
\begin{equation}
\begin{tabular}{ccc}
$L^{\vee }$ & $\overset{F_{\lambda }=\left( \lambda ,\pi \circ q\right) }{
\longrightarrow }$ & $L_{0}^{\vee }\times \Delta $ \\ 
$\downarrow q$ &  & $\downarrow \left( q_{0},id.\right) $ \\ 
$X$ & $\overset{F_{\sigma }=\left( \sigma ,\pi \right) }{\longrightarrow }$
& $X_{0}\times \Delta $ \\ 
$\downarrow \pi $ &  & $\downarrow $ \\ 
$\Delta $ & $=$ & $\Delta $%
\end{tabular}
\label{2.3.1}
\end{equation}
of the deformation $L^{\vee }/X$ of $L_{0}^{\vee }/X_{0}$ with the following
properties:

\begin{enumerate}
\item  The trivialization respects the structure of complex line bundles.
That is, 
\begin{equation}
\left( F_{\lambda }\right) _{*}\left( \chi \right) =\chi  \label{2.3.2}
\end{equation}
where $\chi $ is the Euler vector field for the holomorphic line bundles $%
L^{\vee }$ and $L_{0}^{\vee }\times \Delta $.

\item  For each $x_{0}\in X_{0}$, the restricted map 
\begin{equation*}
F_{\lambda }:\left( q_{0}\circ \lambda \right) ^{-1}\left( x_{0}\right)
\rightarrow \left( q_{0}^{-1}\left( x_{0}\right) \right) \times \Delta
\end{equation*}
is a holomorphic isomorphism.
\end{enumerate}

We prove this ``intuitively obvious'' fact in the Appendix at the end of
this paper.

Our deformation/trivialization $\left( \ref{2.3.1}\right) $ is given by
Kuranishi data 
\begin{equation*}
\xi =\sum\nolimits_{\left| J\right| >0}\xi _{J}t^{J}
\end{equation*}
for which 
\begin{equation*}
\xi _{j}\in A^{0,1}\left( T_{L_{0}^{\vee }}\right)
\end{equation*}
and 
\begin{equation}
L_{\chi }\xi _{j}=L_{\bar{\chi}}\xi _{j}=0.  \label{2.3.3}
\end{equation}

We call a trivialization satisfying $\left( \ref{2.3.1}\right) $-$\left( \ref
{2.3.3}\right) $ a\textit{\ trivialization of line bundles}. We say that the
trivializations $\lambda $ of $L^{\vee }/\Delta $ and $\sigma $ of $X/\Delta 
$ are \textit{compatible} if they make the diagram $\left( \ref{2.3.1}
\right) $ commutative. By an elementary computation in local coordinates,
sections 
\begin{equation*}
\xi _{J}\in A_{L_{0}^{\vee }}^{0,1}\otimes T_{L_{0}^{\vee }}
\end{equation*}
associated to a trivialization of line bundles lie in a subspace 
\begin{equation*}
A\subseteq A_{L_{0}^{\vee }}^{0,1}\otimes T_{L_{0}^{\vee }}
\end{equation*}
comprising the the middle term of an exact sequence 
\begin{equation}
0\rightarrow q_{0}^{-1}\left( A_{X_{0}}^{0,1}\right) \otimes _{\Bbb{C}}\Bbb{C%
}\chi \rightarrow A\rightarrow q_{0}^{-1}\left( A_{X_{0}}^{0,1}\otimes
T_{X_{0}}\right) \rightarrow 0,  \label{2.3.4}
\end{equation}
where $\chi $ is the Euler vector field on $L_{0}^{\vee }$. That is, 
\begin{equation}
A=A_{X_{0}}^{0,1}\left( \frak{D}_{1}\left( L_{0}\right) \right) .
\label{2.3.5}
\end{equation}

Let 
\begin{equation*}
\frak{m}
\end{equation*}
denote the ideal of $0$ in $\Delta $. Notice that the quotient form 
\begin{equation*}
\xi \in A_{X_{0}}^{0,1}\left( \frak{D}_{1}\left( L_{0}\right) \right)
\otimes \frac{\frak{m}}{\frak{m}^{2}}
\end{equation*}
must be $\overline{\partial }$-closed by the integrability conditions in 
\cite{C1}. Its cohomology class in 
\begin{equation*}
H^{1}\left( \frak{D}_{1}\left( L_{0}\right) \right) \otimes \frac{\frak{m}}{%
\frak{m}^{2}}
\end{equation*}
is the first-order deformation of the pair $\left( X_{0},L_{0}\right) $
given by $\left( \ref{2.2}\right) $ (see \cite{AC}). Its symbol is just the
element of $H^{1}\left( T_{X_{0}}\right) \otimes \frac{\frak{m}}{\frak{m}^{2}%
}$ giving the Kodaira-Spencer class for the compatible first-order
deformation of the manifold $X_{0}$.

From \cite{C1} we have:

\begin{lemma}
\label{2.4}A (formal) holomorphic section of $L$ is a power series 
\begin{equation*}
f=\sum\nolimits_{I}t^{I}f_{I}:L_{0}^{\vee }\times \Delta \rightarrow \Bbb{C}
\end{equation*}
with coefficients $f_{I}$ which are $C^{\infty }$-sections of $L_{0}$ such
that 
\begin{equation*}
\left( \bar{\partial }_{L_{0}}-L_{\xi }\right) \left( f\right) =0.
\end{equation*}
\end{lemma}

In what follows, we will work with the bundle 
\begin{equation}
L\oplus \overline{L}  \label{2.4.1}
\end{equation}
and the associated trivialization 
\begin{equation}
F_{\lambda +\overline{\lambda }}:L^{\vee }\oplus \overline{L}^{\vee
}\rightarrow \left( L_{0}^{\vee }\oplus \overline{L}_{0}^{\vee }\right)
\times \Delta .  \label{2.4.2}
\end{equation}
It will be useful to define 
\begin{eqnarray*}
L\oplus \overline{L} &\rightarrow &L\oplus \overline{L} \\
\left( f,g\right) &\rightarrow &\left( \overline{g},\overline{f}\right)
\end{eqnarray*}
and view $\left( \ref{2.4.1}\right) $ as the complexification of real bundle 
\begin{equation*}
L_{\Bbb{R}}=\left\{ f\oplus \overline{f}:f\in L\right\}
\end{equation*}
with almost complex structure given by 
\begin{equation}
\left( f\oplus \overline{f}\mapsto if\oplus \left( -i\right) \overline{f}
\right) .  \label{2.4.3}
\end{equation}

\section{Flat line bundles\label{3}}

Suppose 
\begin{equation*}
X_{0}
\end{equation*}
is a compact K\"{a}hler manifold. We specialize now to deformations of the
trivial bundle 
\begin{equation*}
L_{0}=\mathcal{O}_{X_{0}}.
\end{equation*}
Let 
\begin{equation*}
P\rightarrow Alb\left( X_{0}\right) \times Pic^{0}\left( X_{0}\right)
\end{equation*}
be the Poincare bundle defined as follows. Writing 
\begin{equation*}
\exp \left( a\right) :=e^{2\pi i\cdot a}.
\end{equation*}
and using the the natural inclusion 
\begin{eqnarray}
H_{1}\left( X_{0};\Bbb{Z}\right) &\rightarrow &H^{0}\left( \Omega
_{X_{0}}^{1}\right) ^{\vee }  \label{3.10'''} \\
\varepsilon &\mapsto &\int\nolimits_{\varepsilon }  \notag
\end{eqnarray}
modulo torsion, any $\Bbb{R}$-linear map 
\begin{equation}
\beta :H^{0}\left( \Omega _{X_{0}}^{1}\right) ^{\vee }\rightarrow \Bbb{C}
\label{3.10''}
\end{equation}
gives rise to a quotient space 
\begin{equation}
\frac{\Bbb{C}\times H^{0}\left( \Omega _{X_{0}}^{1}\right) ^{\vee }}{\left\{
\left( z,u+\varepsilon \right) \sim \left( \exp \left( \beta \left(
\varepsilon \right) \right) \cdot z,u\right) \right\} _{\varepsilon \in
H_{1}\left( X_{0};\Bbb{Z}\right) }}  \label{3.10}
\end{equation}
which is a holomorphic line bundle on $Alb\left( X_{0}\right) $. Under the
map 
\begin{equation*}
a:X_{0}\rightarrow \mathrm{Alb}\left( X_{0}\right)
\end{equation*}
given by choosing a basepoint this bundle pulls back to a line bundle $%
L_{\beta }$ and all holomorphic line bundles on $X_{0}$ with trivial first
Chern class are constructed in this way.

Thus we obtain a natural map 
\begin{equation*}
\mathrm{Hom}_{\Bbb{R}}\left( H^{0}\left( \Omega _{X_{0}}^{1}\right) ^{\vee },%
\Bbb{C}\right) \rightarrow Pic^{0}\left( X_{0}\right)
\end{equation*}
which lifts to an $\Bbb{R}$-linear map 
\begin{equation}
\mathrm{Hom}_{\Bbb{R}}\left( H^{0}\left( \Omega _{X_{0}}^{1}\right) ^{\vee },%
\Bbb{C}\right) \rightarrow H^{1}\left( \mathcal{O}_{X_{0}}\right) =:\Delta .
\label{3.10'}
\end{equation}
The kernel of $\left( \ref{3.10'}\right) $ is the space $H^{0}\left( \Omega
_{X_{0}}^{1}\right) $, that is, the set of $\Bbb{C}$-linear maps 
\begin{equation*}
H^{0}\left( \Omega _{X_{0}}^{1}\right) ^{\vee }\rightarrow \Bbb{C}.
\end{equation*}
Using the inclusion $\left( \ref{3.10'''}\right) $ we identify 
\begin{equation*}
\begin{array}{c}
\mathrm{Hom}_{\Bbb{R}}\left( H^{0}\left( \Omega _{X_{0}}^{1}\right) ^{\vee },%
\Bbb{C}\right) \cong H^{1}\left( X_{0};\Bbb{C}\right) \\ 
\beta \mapsto d\beta
\end{array}
\end{equation*}
so that the mapping $\left( \ref{3.10'}\right) $ becomes the natural
surjection 
\begin{equation}
H^{1}\left( X_{0};\Bbb{C}\right) \rightarrow H^{1}\left( \mathcal{O}%
_{X_{0}}\right)  \label{natsurj}
\end{equation}
with kernel 
\begin{equation*}
\mathrm{Hom}_{\Bbb{C}}\left( H^{0}\left( \Omega _{X_{0}}^{1}\right) ^{\vee },%
\Bbb{C}\right) =H^{0}\left( \Omega _{X_{0}}^{1}\right) .
\end{equation*}

As Deligne has pointed out (see \cite{S2}, page 367), the space 
\begin{eqnarray*}
M &=&H^{1}\left( X_{0};\Bbb{C}^{*}\right) \\
&=&Hom\left( \pi _{1}\left( X_{0}\right) ,\Bbb{C}^{*}\right)
\end{eqnarray*}
has the structure of a quaternionic vector space, where the real tangent
space 
\begin{equation*}
T_{M}\left( \Bbb{R}\right) =M\times H^{1}\left( X_{0};\Bbb{C}\right)
\end{equation*}
and:

\begin{enumerate}
\item  Multiplication of the coefficients in $H^{1}\left( X_{0};\Bbb{C}%
\right) $ by $\sqrt{-1}$ induces a transformation $k$ of the $\Bbb{R}$
-vector space 
\begin{equation*}
H^{1}\left( X_{0};\Bbb{C}\right)
\end{equation*}
giving the (usual) complex structure on $M$.

\item  On the other hand we identify the complex vector space $H^{1}\left( 
\mathcal{O}_{X_{0}}\right) $ with the subspace 
\begin{equation*}
H^{1}\left( X_{0};\Bbb{R}\right) \subseteq H^{1}\left( X_{0};\Bbb{C}\right)
\end{equation*}
via $\left( \ref{natsurj}\right) $ and identify the complex vector space $%
H^{0}\left( \Omega _{X_{0}}^{1}\right) $ with 
\begin{equation*}
H^{1}\left( X_{0};i\cdot \Bbb{R}\right) \subseteq H^{1}\left( X_{0};\Bbb{C}%
\right)
\end{equation*}
via the map 
\begin{equation*}
H^{0}\left( \Omega _{X_{0}}^{1}\right) \rightarrow \frac{H^{1}\left( X_{0};%
\Bbb{C}\right) }{H^{1}\left( X_{0};\Bbb{R}\right) }=H^{1}\left( X_{0};\sqrt{%
-1}\cdot \Bbb{R}\right) .
\end{equation*}
Multiplication by $\sqrt{-1}$ in each of these complex vector spaces give a
complex structure $j$ on $H^{1}\left( X_{0};\Bbb{C}\right) $ which respects
the direct-sum decomposition 
\begin{equation*}
H^{1}\left( X_{0};\Bbb{C}\right) =H^{1}\left( X_{0};\Bbb{R}\right) \oplus
H^{1}\left( X_{0};i\cdot \Bbb{R}\right) ,
\end{equation*}
that is, the decomposition 
\begin{eqnarray*}
H^{1}\left( X_{0};\Bbb{C}^{*}\right) =\mathrm{Hom}_{\Bbb{R}}\left(
H_{1}\left( X_{0};\Bbb{Z}\right) ,\Bbb{U}\left( 1\right) \right) \times 
\mathrm{Hom}_{\Bbb{R}}\left( H^{1}\left( X_{0};\Bbb{Z}\right) ,\Bbb{R}%
_{+}^{*}\right) \\
=\exp \left( H^{1}\left( X_{0};\Bbb{R}\right) \right) \times \exp \left(
H^{1}\left( X_{0};\sqrt{-1}\cdot \Bbb{R}\right) \right)
\end{eqnarray*}
\end{enumerate}

\noindent Later we use 
\begin{equation*}
\beta \left( u,\zeta \right)
\end{equation*}
to denote the dependence of $\beta $ on the holomorphic coordinate $\zeta $
of $\Delta =H^{1}\left( \mathcal{O}_{X_{0}}\right) $ and the holomorphic
coordinate $u$ of $H^{0}\left( \Omega _{X_{0}}^{1}\right) $. To relate the
two complex structures 
\begin{equation*}
\begin{array}{ccc}
& H^{1}\left( X_{0};\Bbb{C}\right) &  \\ 
& = &  \\ 
H^{0}\left( \Omega _{X_{0}}^{1}\right) & \oplus & H^{1}\left( \mathcal{O}%
_{X_{0}}\right) \\ 
\downarrow &  & \downarrow \\ 
H^{1}\left( X_{0};\sqrt{-1}\cdot \Bbb{R}\right) & \oplus & H^{1}\left( X_{0};%
\Bbb{R}\right) \\ 
& = &  \\ 
& H^{1}\left( X_{0};\Bbb{C}\right) & 
\end{array}
\end{equation*}
where the vertical maps are given by 
\begin{equation*}
\begin{array}{ccc}
\omega & + & \eta \\ 
\downarrow &  & \downarrow \\ 
\frac{\omega -\overline{\omega }}{\sqrt{2}} & + & \frac{\eta +\overline{\eta 
}}{\sqrt{2}}
\end{array}
.
\end{equation*}
So 
\begin{eqnarray*}
k\cdot j\cdot \left( \omega +\eta \right) &=&k\cdot \left( \sqrt{-1}\cdot
\omega +\sqrt{-1}\cdot \eta \right) \\
&=&k\cdot \left( \sqrt{-1}\cdot \frac{\omega +\overline{\omega }}{\sqrt{2}}+%
\sqrt{-1}\cdot \frac{\eta -\overline{\eta }}{\sqrt{2}}\right) \\
&=&-\left( \frac{\eta -\overline{\eta }}{\sqrt{2}}+\frac{\omega +\overline{%
\omega }}{\sqrt{2}}\right) \\
&=&\left( \overline{\eta }+\left( -\overline{\omega }\right) \right)
\end{eqnarray*}
and 
\begin{eqnarray*}
j\cdot k\cdot \left( \omega +\eta \right) &=&j\cdot k\cdot \left( \frac{%
\omega -\overline{\omega }}{\sqrt{2}}+\frac{\eta +\overline{\eta }}{\sqrt{2}}%
\right) \\
&=&j\cdot \left( \sqrt{-1}\cdot \frac{\eta +\overline{\eta }}{\sqrt{2}}+%
\sqrt{-1}\cdot \frac{\omega -\overline{\omega }}{\sqrt{2}}\right) \\
&=&j\cdot \left( \frac{\sqrt{-1}\cdot \eta -\overline{\sqrt{-1}\cdot \eta }}{%
\sqrt{2}}+\frac{\sqrt{-1}\cdot \omega +\overline{\sqrt{-1}\cdot \omega }}{%
\sqrt{2}}\right) \\
&=&j\cdot \left( -\overline{\sqrt{-1}\cdot \eta }+\overline{\sqrt{-1}\cdot
\omega }\right) \\
&=&\left( \left( -\overline{\eta }\right) +\overline{\omega }\right) .
\end{eqnarray*}
So $i=jk$, $j$ and $k$ give the real vector space $H^{1}\left( X_{0};\Bbb{C}%
\right) $ the structure of a quaternionic vector space.

Via the pullback of the Poincare bundle 
\begin{equation*}
L=a^{*}P
\end{equation*}
we obtain a distinguished $C^{\infty }$-trivialization 
\begin{eqnarray}
L &\rightarrow &\left( \Bbb{C}\times X_{0}\right) \times H^{1}\left( X_{0};%
\Bbb{C}\right) .  \label{3.0} \\
\left( z,\left( u,\zeta \right) \right) &\mapsto &\left( \exp \left( -\beta
\right) \cdot z,\left( u,\zeta \right) \right)  \notag
\end{eqnarray}
Thus for this deformation/trivialization we have the Kuranishi data over $%
X_{0}$ given by 
\begin{equation}
\xi =\overline{\partial }\beta \cdot \chi  \label{3.1}
\end{equation}
where, abusing notation, we denote $\beta \circ a$ simply as $\beta $.

Since the complex structure varies only in the directions $H^{1}\left( 
\mathcal{O}_{X_{0}}\right) \subseteq H^{1}\left( X_{0};\Bbb{C}\right) $ we
will let 
\begin{equation*}
\Delta =H^{1}\left( \mathcal{O}_{X_{0}}\right)
\end{equation*}
and, abusing notation, write 
\begin{equation*}
L/\left( X_{0}\times \Delta \right) .
\end{equation*}
However we will usually wish to choose the $\beta $ from another subspace of 
$H^{1}\left( X_{0};\Bbb{C}\right) $, namely 
\begin{equation*}
H^{1}\left( X_{0};\Bbb{R}\right)
\end{equation*}
These two subspaces are identified, of course, under the restriction of the
natural projection 
\begin{equation*}
H^{1}\left( X_{0};\Bbb{C}\right) \rightarrow H^{1}\left( \mathcal{O}
_{X_{0}}\right)
\end{equation*}
to $H^{1}\left( X_{0};\Bbb{R}\right) .$

$\left( \ref{3.0}\right) $ induces the isomorphism 
\begin{equation}
\begin{array}{c}
A_{X_{0}\times \Delta /\Delta }^{*}\left( L\oplus \overline{L}\right) \\ 
\downarrow \tilde{\varphi} \\ 
A_{X_{0}}^{*}\left( \mathcal{O}_{X_{0}}\oplus \overline{\mathcal{O}_{X_{0}}}
\right) \otimes A_{\Delta }^{0}
\end{array}
\label{3.e}
\end{equation}
and therefore the correspondence of operators 
\begin{eqnarray}
\overline{\partial }_{L} &\leftrightarrow &\overline{\partial }-\overline{
\partial }\beta \wedge  \label{3.2} \\
\partial _{\overline{L}} &\leftrightarrow &\partial -\partial \overline{
\beta }\wedge .  \notag
\end{eqnarray}
(Notice that, since $\beta $ is linear on $\mathrm{Alb}\left( X_{0}\right) $
, $\partial \overline{\partial }\beta =0$.)

Let 
\begin{equation*}
\pi :X_{0}\times \mathrm{Pic}^{0}X_{0}\rightarrow \mathrm{Pic}^{0}X_{0}
\end{equation*}
be the standard projection and 
\begin{equation*}
\rho :H^{1}\left( \mathcal{O}_{X_{0}}\right) \rightarrow \mathrm{Pic}%
^{0}X_{0}
\end{equation*}
the standard quotient map. Furthermore just as as in Proposition 2.3 of \cite
{GL2} we have from $\left( \ref{3.2}\right) $ that 
\begin{equation}
R^{q}\pi _{*}\left( \Omega _{X_{0}}^{p}\otimes a^{*}P\right) \cong \mathcal{H%
}^{q}\left( A_{X_{0}}^{p,*}\otimes \mathcal{O}_{\Delta },\overline{\partial }%
-\overline{\partial }\beta \wedge \right)  \label{3.3}
\end{equation}
where 
\begin{equation*}
\pi :X_{0}\times \Delta \rightarrow \Delta
\end{equation*}
is the standard projection.

\section{The line-bundle metric\label{4}}

For any flat line bundle $L_{0}$%
\begin{equation}
c_{1}\left( L_{0}\right) =\frac{i}{2\pi }\partial _{X_{0}}\overline{\partial 
}_{X_{0}}\log \kappa _{0}=0  \label{4.1}
\end{equation}
where 
\begin{equation*}
\kappa _{0}\left( f,g\right)
\end{equation*}
is the flat hermitian metric on $L_{0}$. Notice that $\kappa _{0}$ is given
by the complexification 
\begin{equation}
\left\langle \ ,\ \right\rangle _{0}  \label{4.2}
\end{equation}
of a real inner product on 
\begin{equation*}
L_{\Bbb{R},0}\times X_{0}
\end{equation*}
(for which $\left( \ref{2.4.3}\right) $ is an isometry) via the rule 
\begin{equation*}
\kappa _{0}\left( f,g\right) =\left\langle f\oplus \overline{f},g\oplus 
\overline{g}\right\rangle _{0}.
\end{equation*}

For any deformation/trivialization 
\begin{equation*}
\left( \lambda ,\pi \circ \left( q\oplus \overline{q}\right) \right)
:L\oplus \overline{L}\rightarrow \left( L_{0}\oplus \overline{L_{0}}\right)
\times \Delta ,
\end{equation*}
for example for the one constructed in \S \ref{3}, define a real inner
product on $L_{\Bbb{R}}\subseteq L\oplus \overline{L}$ by the rule 
\begin{equation}
\left\langle \ ,\ \right\rangle =\lambda ^{*}\left\langle \ ,\ \right\rangle
_{0}  \label{4.11}
\end{equation}
and extend $\kappa _{0}$ to a hermitian metric $\kappa $ on $L\oplus 
\overline{L}$ in the standard way. This hermitian metric is given on $L$ by
the rule that 
\begin{equation*}
\kappa \left( f,g\right) =\left\langle f,\overline{g}\right\rangle .
\end{equation*}

Thus if 
\begin{equation*}
d=\dim _{\Bbb{C}}X_{0}
\end{equation*}
the (real) ``star'' operator 
\begin{equation}
\#:A^{*}\left( L\oplus \overline{L}\right) \rightarrow A^{2d-*}\left(
L\oplus \overline{L}\right)  \label{4.12}
\end{equation}
corresponds under the isomorphism $\left( \ref{3.e}\right) $ to the
pull-back of the real star operator 
\begin{equation*}
\#_{0}:A^{*}\left( L_{0}\oplus \overline{L_{0}}\right) \rightarrow
A^{2d-*}\left( L_{0}\oplus \overline{L_{0}}\right) .
\end{equation*}

Let 
\begin{equation*}
D_{0}:A^{*}\left( L_{0}\oplus \overline{L_{0}}\right) \rightarrow
A^{*+1}\left( L_{0}\oplus \overline{L_{0}}\right)
\end{equation*}
be the distinguished real connection induced by $\kappa _{0}$ and the
K\"{a}hler metric on $X_{0}$. This connection respects the direct sum
decomposition $L_{0}\oplus \overline{L_{0}}$). Then the (relative)
connection 
\begin{equation*}
D_{L\oplus \overline{L}}:A_{X_{0}\times \Delta /\Delta }^{*}\left( L\oplus 
\overline{L}\right) \rightarrow A_{X_{0}\times \Delta /\Delta }^{*+1}\left(
L\oplus \overline{L}\right)
\end{equation*}
induced by $\kappa $ and the ``constant'' relative K\"{a}hler metric on $%
X_{0}\times \Delta /\Delta $ corresponds under the isomorphism $\tilde{%
\varphi}$ given in $\left( \ref{2.4.2}\right) $ to a (relative) connection 
\begin{equation*}
\begin{array}{c}
\left( A^{*}\left( L_{0}\oplus \overline{L_{0}}\right) \otimes A_{\Delta
}^{0}\right) \\ 
\downarrow D \\ 
\left( A^{*+1}\left( L_{0}\oplus \overline{L_{0}}\right) \otimes A_{\Delta
}^{0}\right) .
\end{array}
\end{equation*}
So for the deformation in \S \ref{3} we have 
\begin{equation*}
D=d-\left( \left( d\beta \wedge \right) \oplus \left( d\overline{\beta }%
\wedge \right) \right) .
\end{equation*}

Also $L\oplus \overline{L}$ is the complexification of a real bundle $L_{%
\Bbb{R}}$ spanned by 
\begin{eqnarray*}
Re\alpha &=&\frac{\alpha \oplus \overline{\alpha }}{2} \\
Im\alpha &=&\frac{i\alpha \oplus \overline{i\alpha }}{2}
\end{eqnarray*}
for sections $\alpha $ of $L$. $D_{L\oplus \overline{L}}$ is the
complexification of a real connection on $L_{\Bbb{R}}$ given under the
isomorphism $\tilde{\varphi}$ of $\left( \ref{3.e}\right) $ by the formula 
\begin{equation}
d\left( Re\alpha _{0},\ Im\alpha _{0}\right) +\left( Re\alpha _{0},\
Im\alpha _{0}\right) \wedge \left( 
\begin{array}{cc}
Re\left( d\beta \right) & -Im\left( d\beta \right) \\ 
Im\left( d\beta \right) & Re\left( d\beta \right)
\end{array}
\right) .  \label{4.6}
\end{equation}
(Here $Re\left( d\beta \right) $ and $Im\left( d\beta \right) $ are
respectively the real and imaginary parts of the one-form $d\beta $.)

\section{Harmonic forms\label{5}}

To proceed further with our analysis of harmonic forms, we continue to
restrict our considerations of \S \ref{4} to the case 
\begin{equation*}
L_{0}=\mathcal{O}_{X_{0}}
\end{equation*}
with the standard flat metric $\kappa _{0}$. Thus the connection $\left( \ref
{4.6}\right) $ gives for each 
\begin{equation*}
\beta \in \mathrm{Hom}_{\Bbb{R}}\left( H^{0}\left( \Omega
_{X_{0}}^{1}\right) ^{\vee },\Bbb{C}\right) =H^{1}\left( X_{0};\Bbb{C}\right)
\end{equation*}
a flat connection on the trivial bundle 
\begin{equation*}
X_{0}\times \Bbb{R}^{2}\cong L_{\Bbb{R},\beta }\subseteq L_{\beta }\oplus 
\overline{L_{\beta }}.
\end{equation*}
Following \cite{S1} our goal is to give explicitly the identification
between local systems $L_{\Bbb{R},\beta }$ and Higgs bundle structures
implicit in the choice of $\beta $.

Continuing with 
\begin{equation*}
L_{0}=\mathcal{O}_{X_{0}}
\end{equation*}
the pasting functions defining $L_{\beta }$ are locally constant so that 
\begin{equation*}
\partial _{\beta }:A_{X_{0}}^{p,q}\left( L_{\beta }\right) \rightarrow
A_{X_{0}}^{p+1,q}\left( L_{\beta }\right)
\end{equation*}
is well defined. Also, under the isomorphism 
\begin{equation}
\tilde{\varphi}:A_{X_{0}}^{*}\left( L_{\beta }\right) \oplus
A_{X_{0}}^{*}\left( \overline{L_{\beta }}\right) \leftrightarrow
A_{X_{0}}^{*}\oplus A_{X_{0}}^{*}  \label{5.4}
\end{equation}
given in $\left( \ref{3.e}\right) $ we can complete the correspondence 
\begin{equation*}
\overline{\partial }_{\beta }\leftrightarrow \overline{\partial }-\overline{
\partial }\beta \wedge .
\end{equation*}
given as in $\left( \ref{3.2}\right) $ as follows: 
\begin{equation*}
d_{\beta }\oplus d_{\overline{\beta }}\leftrightarrow \left( d-d\beta \wedge
\right) \oplus \left( d-d\overline{\beta }\wedge \right) =D
\end{equation*}
where 
\begin{eqnarray*}
d_{\beta } &=&\partial _{\beta }+\overline{\partial }_{\beta } \\
d_{\overline{\beta }} &=&\partial _{\overline{\beta }}+\overline{\partial }_{%
\overline{\beta }}
\end{eqnarray*}
with 
\begin{eqnarray*}
\partial _{\beta } &\leftrightarrow &\partial -\partial \beta \wedge \\
\partial _{\overline{\beta }} &\leftrightarrow &\partial -\partial \overline{
\beta }\wedge \\
\overline{\partial }_{\beta } &\leftrightarrow &\overline{\partial }-%
\overline{\partial }\beta \wedge \\
\overline{\partial }_{\overline{\beta }} &\leftrightarrow &\overline{
\partial }-\overline{\partial }\overline{\beta }\wedge .
\end{eqnarray*}

From $\left( \ref{4.6}\right) $ we have $D\left( Re\alpha ,\ Im\alpha
\right) $ given by the formula 
\begin{equation*}
\left( dRe\alpha ,\ dIm\alpha \right) +\left( Re\alpha ,\ Im\alpha \right)
\wedge \left( 
\begin{array}{cc}
Re\left( d\beta \right) & -Im\left( d\beta \right) \\ 
Im\left( d\beta \right) & Re\left( d\beta \right)
\end{array}
\right) .
\end{equation*}
Let 
\begin{equation*}
D=D^{\prime }+D^{\prime \prime }
\end{equation*}
be the decomposition of $D$ into its $\left( 1,0\right) $- and $\left(
0,1\right) $-components respectively. We have 
\begin{eqnarray*}
D^{\prime } &=&\left( \partial -\partial \beta \wedge \right) \oplus \left(
\partial -\partial \overline{\beta }\wedge \right) \\
D^{\prime \prime } &=&\left( \overline{\partial }-\overline{\partial }\beta
\wedge \right) \oplus \left( \overline{\partial }-\overline{\partial }%
\overline{\beta }\wedge \right) .
\end{eqnarray*}

Also for the flat metric connection 
\begin{equation*}
d_{\beta }\oplus d_{\overline{\beta }}:A^{0}\left( L_{\beta }\oplus 
\overline{L_{\beta }}\right) \rightarrow A^{1}\left( L_{\beta }\oplus 
\overline{L_{\beta }}\right)
\end{equation*}
and the $\Bbb{C}$-linear isomorphism 
\begin{equation*}
\#:L_{\beta }\rightarrow \overline{L_{\beta }}^{\vee }
\end{equation*}
induced by the metric as in $\left( \ref{4.12}\right) $ gives a
distinguished real $C^{\infty }$-isomorphism 
\begin{equation}
\#:A^{p,q}\left( L_{\beta }\oplus \overline{L_{\beta }}\right) \rightarrow
A^{n-q,n-p}\left( \overline{L_{\beta }}^{\vee }\oplus L_{\beta }^{\vee
}\right)  \label{4.4}
\end{equation}
corresponding under $\tilde{\varphi}$ in $\left( \ref{5.4}\right) $ to 
\begin{equation*}
\begin{array}{c}
\ast :A^{p,q}\oplus A^{p,q}\rightarrow A^{n-q,n-p}\oplus A^{n-q,n-p}. \\ 
\alpha \oplus \beta \mapsto *\beta \oplus *\alpha
\end{array}
\end{equation*}
(Caution: The definition of the star operator used in this paper is the
complexification of the real star operator on real differential forms
induced by the Riemannian metric on $X_{0}$. In particular, this operator is 
$\Bbb{C}$-linear.)

Thus 
\begin{equation*}
D^{*}=-*\circ D\circ *.
\end{equation*}

On the other hand, with respect to a local flat framing of $L_{\beta }\oplus 
\overline{L_{\beta }}$ one has that $\#$ is just given by the usual star
operator on the K\"{a}hler manifold $X_{0}$ and the connection $D_{L_{\beta
}\oplus \overline{L_{\beta }}}$, on $A^{*}\left( L_{\beta }\oplus \overline{
L_{\beta }}\right) $ is given by the usual exterior derivative. Thus $\Delta
_{D_{L_{\beta }\oplus \overline{L_{\beta }}}}$ preserves type and so, for
example 
\begin{equation}
\partial _{\beta }\circ \overline{\partial }_{\beta }^{*}+\overline{\partial 
}_{\beta }^{*}\circ \partial _{\beta }=0.  \label{5.0'}
\end{equation}
Therefore the Laplacian 
\begin{equation*}
\Delta _{D}:\left( A_{X_{0}}^{*}\oplus A_{X_{0}}^{*}\right) \rightarrow
\left( A_{X_{0}}^{*}\oplus A_{X_{0}}^{*}\right)
\end{equation*}
corresponding to $\Delta _{D_{L_{\beta }\oplus \overline{L_{\beta }}}}$
under $\tilde{\varphi}$ also preserves type

Using $\left( \ref{4.4}\right) $ we have by an easy calculation that 
\begin{eqnarray*}
\left( d_{L}\right) ^{*} &=&\leftrightarrow d^{*}-\left( d\beta \wedge
\right) ^{*}=-*\left( d-d\overline{\beta }\wedge \right) * \\
\left( d_{\overline{L}}\right) ^{*} &\leftrightarrow &-*\left( d-d\beta
\wedge \right) *.
\end{eqnarray*}
Thus 
\begin{eqnarray}
\left( D^{\prime }\right) ^{*} &=&\left( \partial ^{*}-*\left( \overline{%
\partial }\overline{\beta }\wedge \right) *\right) \oplus \left( \partial
^{*}-*\left( \overline{\partial }\beta \wedge \right) *\right)  \label{5.1'}
\\
\left( D^{\prime \prime }\right) ^{*} &=&\left( \overline{\partial }
^{*}-*\left( \partial \overline{\beta }\wedge \right) *\right) \oplus \left( 
\overline{\partial }^{*}-*\left( \partial \beta \wedge \right) *\right) . 
\notag
\end{eqnarray}
Since $\Delta _{D}$ preserves type, keeping track of summands we have 
\begin{equation*}
\Delta _{D}=\Delta _{D}^{1}\oplus \Delta _{D}^{2}=\left( \Delta _{D^{\prime
}}^{1}+\Delta _{D^{\prime \prime }}^{1}\right) \oplus \left( \Delta
_{D^{\prime }}^{2}+\Delta _{D^{\prime \prime }}^{2}\right) .
\end{equation*}
Also 
\begin{eqnarray*}
\overline{\Delta _{D^{\prime }}^{1}} &=&\Delta _{D^{\prime \prime }}^{2} \\
\overline{\Delta _{D^{\prime \prime }}^{1}} &=&\Delta _{D^{\prime }}^{2}.
\end{eqnarray*}
But, since the metric $\kappa $ is flat, the operator $\Delta _{D^{\prime
\prime }}^{1}$ is the restriction to the first summand of a real operator
whose restriction to the second summand is $\Delta _{D^{\prime \prime }}^{2}$
. So 
\begin{equation}
\Delta _{D}^{1}=2\Delta _{D^{\prime }}^{1}=2\Delta _{D^{\prime \prime
}}^{1},\ \Delta _{D}^{2}=2\Delta _{D^{\prime }}^{2}=2\Delta _{D^{\prime
\prime }}^{2}.  \label{5del}
\end{equation}

Thus:

\begin{theorem}
\label{p5}i) The harmonic forms in $A_{X_{0}}^{*}\left( L_{\beta }\right) $
correspond under $\left( \ref{3.3}\right) $ to forms $\alpha \in
A_{X_{0}}^{*}$ for which 
\begin{eqnarray}
\overline{\partial }\alpha =\overline{\partial }\beta \left( \zeta \right)
\wedge \alpha  \label{5.5} \\
\partial *\alpha =\partial \overline{\beta }\left( \zeta \right) \wedge
\left( *\alpha \right) .  \notag
\end{eqnarray}
Also 
\begin{equation*}
\overline{H^{q}\left( \Omega _{X_{0}}^{p}\otimes L_{\beta }\right) }\cong
H^{p}\left( \Omega _{X_{0}}^{q}\otimes L_{\overline{\beta }}\right)
\end{equation*}
and 
\begin{equation*}
\sum\nolimits_{p+q=r}H^{q}\left( \Omega _{X_{0}}^{p}\otimes L_{\beta
}\right) \oplus \overline{H^{q}\left( \Omega _{X_{0}}^{p}\otimes L_{\beta
}\right) }\cong H^{r}\left( X_{0};L_{\Bbb{R},\beta }\right) \otimes \Bbb{C}
\end{equation*}
where the right-hand group is the cohomology of $X_{0}$ with coefficients in
the real local system $L_{\Bbb{R}}$ corresponding to the flat connection $%
\left( \ref{4.6}\right) $. (See \cite{S1}, Lemma 1.2, and \cite{S2}, \S 5.)
\end{theorem}

ii) The $\partial \overline{\partial }$-lemma holds for the operators $%
\partial _{\beta }$ and $\overline{\partial }_{\beta }$, that is, if $\alpha
\in A^{p,q}\left( L_{\beta }\right) $ is $\partial _{\beta }$-closed (exact)
and $\overline{\partial }_{\beta }$-exact (closed), then 
\begin{equation*}
\alpha \in \partial _{\beta }\overline{\partial }_{\beta }A^{p-1,q-1}\left(
L_{\beta }\right)
\end{equation*}

\begin{proof}
i) By $\left( \ref{5del}\right) $ $D^{\prime }$-cohomology and $D^{\prime
\prime }$-cohomology coincide.

ii) The usual proof of the $\partial \overline{\partial }$-lemma goes
through since we have the Green's operator, harmonic decomposition, etc. as
usual and $\left( \ref{5.0'}\right) $
\end{proof}

There are two essential points here. The first is that the moduli space of
local systems $L_{\Bbb{R}}$ is simply the moduli space of representations of
the fundamental group 
\begin{equation*}
\rho :\pi _{1}\left( X_{0}\right) \rightarrow GL\left( 2,\Bbb{R}\right)
\end{equation*}
and so is independent of the complex structure on $X_{0}$. The
representations constructed above consist in all those which factor through
elements of 
\begin{equation}
\beta \in \mathrm{Hom}_{\Bbb{R}}\left( H_{1}\left( X_{0};\Bbb{Z}\right) ,%
\Bbb{C}\right) =H^{1}\left( X_{0};\Bbb{C}\right) \overset{\exp }{
\longrightarrow }H^{1}\left( X_{0};\Bbb{C}^{*}\right) \subseteq H^{1}\left(
X_{0};GL\left( 2,\Bbb{R}\right) \right) .  \label{beta}
\end{equation}
So suppose that we have a smooth deformation 
\begin{equation*}
X_{t}
\end{equation*}
of $X_{0}$. Under the isomorphism 
\begin{equation*}
H^{1}\left( X_{t};\Bbb{C}\right) =H^{1}\left( X_{0};\Bbb{C}\right)
\end{equation*}
induced by the Gauss-Manin connection, we fix $\beta $ as in $\left( \ref
{beta}\right) $. Then the dimension of the vector spaces 
\begin{equation*}
H^{i}\left( X_{t};L_{\Bbb{R},\beta }\right)
\end{equation*}
is locally constant in $t$. So if $\beta \in H^{1}\left( X_{0};\Bbb{R}
\right) $ then by Theorem \ref{p5} and semi-continuity, the dimension of
each of the summands 
\begin{equation*}
H^{q}\left( X_{t};\Omega _{X_{0}}^{p}\otimes L_{\beta }\right)
\end{equation*}
of 
\begin{equation*}
H^{p+q}\left( X_{t};L_{\Bbb{R},\beta }\right)
\end{equation*}
is locally constant in $t$.

The second essential point is the interplay between the two complex
structures $j$ and $k$ defined in \S \ref{3}$.$ First we think of the local
systems $L_{\Bbb{R},\beta }$ as parametrized by the complex manifold 
\begin{equation*}
H^{1}\left( X_{0};\Bbb{C}^{*}\right)
\end{equation*}
with complex structure $k$. Taking, for example, a fixed triangulation of
the topological manifold $X_{0}$, we fix the complex of $\Bbb{C}^{*}$-valued
cochains and write the family of coboundary maps operators as
operator-valued functions on $H^{1}\left( X_{0};\Bbb{C}^{*}\right) $. Since
the loci of constant rank of these operators give a stratification of $%
H^{1}\left( X_{0};\Bbb{C}^{*}\right) $ in which each stratum is a complex
submanifold with respect to the complex structure $k$, each locus 
\begin{equation*}
h^{r}\left( L_{\Bbb{R},\beta }\right) \geq m_{r}
\end{equation*}
is stratified by analytic subvarieties of $H^{1}\left( X_{0};\Bbb{C}\right) $
with respect to the complex structure $k$.

On the other hand, the collection of $\beta $ such that the rank of the
Dolbeault cohomology 
\begin{equation}
h^{q}\left( \Omega _{X_{0}}^{p}\otimes L_{\beta }\right) \geq m_{p,q}
\label{geq}
\end{equation}
is an analytic subvariety with respect to the complex structure $j$ given in
\S \ref{3} since $h^{q}\left( \Omega _{X_{0}}^{p}\otimes L_{\beta }\right) $
is constant on fibers of the map 
\begin{equation*}
H^{1}\left( X_{0};\Bbb{C}\right) \rightarrow H^{1}\left( \mathcal{O}%
_{X_{0}}\right)
\end{equation*}
and $-j$ is compatible with the usual complex structure on $H^{1}\left( 
\mathcal{O}_{X_{0}}\right) $.

By semi-continuity, wherever $h^{p+q}\left( L_{\Bbb{R},\beta }\right) $ is
locally constant, each $h^{q}\left( \Omega _{X_{0}}^{p}\otimes L_{\beta
}\right) $ is too, and so the choice of 
\begin{equation*}
m_{r}=\sum\nolimits_{p+q=r}m_{p,q}
\end{equation*}
locally determines each $m_{p,q}$. Thus the locus at which $h^{p+q}\left( L_{%
\Bbb{R},\beta }\right) \geq m_{p+q}$ is also analytic with respect to the
complex structure $j$. So (following Simpson and Deligne in \cite{S2}, page
367) any irreducible component of this locus is a translate of an $\Bbb{R}$
-linear subspace of $H^{1}\left( X_{0};\Bbb{C}\right) $ which is invariant
under both $j$ and $k$. (Repeating Deligne's proof from \cite{S2}, at a
general point, any irreducible component is locally the graph of a
quaternionic-holomorphic function on a quaterionic manifold. So the Hessian
matrix $Q$ of the function, considered as a quadratic form on the
quaternionic tangent space satisfies 
\begin{eqnarray*}
jk\cdot Q\left( u,v\right) &=&j\cdot Q\left( u,k\cdot v\right) =Q\left(
j\cdot u,k\cdot v\right) \\
&=&k\cdot Q\left( j\cdot u,v\right) =k\cdot j\cdot Q\left( u,v\right)
=-jk\cdot Q\left( u,v\right)
\end{eqnarray*}
and so $Q=0$.)

Since each locus 
\begin{equation*}
h^{q}\left( \Omega _{X_{0}}^{p}\otimes L_{\beta }\right) \geq m_{p,q}
\end{equation*}
is (locally) independent of the complex structure, we can choose a nearby $%
X_{t}$ defined over $\overline{\Bbb{Q}}$. Thus, by Theorem 3.3 of \cite{S2}$%
: $

\begin{theorem}
\label{p5'}For any positive integer $m_{p,q}$, the set of all $L\in \mathrm{%
\ \ Pic}^{0}X_{0}$ for which 
\begin{equation}
h^{q}\left( \Omega _{X_{0}}^{p}\otimes L_{\beta }\right) \geq m_{p,q}
\label{lb}
\end{equation}
is a union of torsion translates of abelian subvarieties of $\mathrm{Pic}%
^{0}X_{0}$. If $X_{0}$ deforms to a smooth $X_{t}$ such that 
\begin{equation*}
\mathrm{Pic}^{0}\left( X_{t}\right)
\end{equation*}
is a simple abelian variety, then the only subtori of $\mathrm{Pic}%
^{0}\left( X_{0}\right) $ which are given by $\left( \ref{lb}\right) $ are
trivial subtori, namely the zero-dimensional subtorus and $\mathrm{Pic}%
^{0}\left( X_{0}\right) $ itself.
\end{theorem}

We next explore the local scheme structure defined by the condition $\left( 
\ref{lb}\right) $ in $\mathrm{Pic}^{0}\left( X_{0}\right) $.

\section{The stalk of the relative cohomology sheaf\label{6}}

We now let $\zeta $ be denote a complex linear coordinate system on $%
H^{1}\left( \mathcal{O}_{X_{0}}\right) .$ We now wish to use $\left( \zeta ,%
\overline{\zeta }\right) $ as formal parameters for the deformation of a
line bundle $\left\{ L_{0}\right\} \in \mathrm{Pic}^{0}\left( X_{0}\right) $
given by 
\begin{equation*}
a^{*}P
\end{equation*}
where as before $P$ is the Poincare bundle on $\mathrm{Alb}X_{0},$ that is, 
\begin{equation*}
\zeta =0
\end{equation*}
is the point $\left\{ L_{0}\right\} \in \mathrm{Pic}^{0}\left( X_{0}\right) $
. As above we compute with the coboundary operators given by the Kuranishi
data and the sheaf of relative differentials. As in \cite{C1} we will
consider this as a complex over the complete local ring 
\begin{equation*}
\Bbb{C}\left[ \left[ \zeta \right] \right] =\frac{\Bbb{C}\left[ \left[ \zeta
,\overline{\zeta }\right] \right] }{\left\{ \overline{\zeta }\right\} }.
\end{equation*}

Let $\frak{m}$ denote the maximal ideal at $0$ in $Spec\Bbb{C}\left[ \left[
\zeta \right] \right] $ which corresponds to the point $\left\{
L_{0}\right\} \in \mathrm{Pic}^{0}\left( X_{0}\right) $. We can write the
graded ring 
\begin{equation}
A^{p,q}\left( X_{0}\right) \otimes \Bbb{C}\left[ \left[ \zeta \right]
\right] =\bigoplus\nolimits_{m=0}^{\infty }A^{p,q}\left( X_{0}\right)
\otimes \frac{\frak{m}^{m}}{\frak{m}^{m+1}}.  \label{6.0}
\end{equation}
So, if $\frak{A}$ is a homogeneous ideal, we have an induced grading on 
\begin{equation*}
A^{p,q}\left( X_{0}\right) \otimes \frac{\Bbb{C}\left[ \left[ \zeta \right]
\right] }{\frak{A}}
\end{equation*}
The coboundary operator 
\begin{equation*}
\overline{\partial }_{\beta }:A^{p,q}\left( X_{0}\right) \otimes \frac{\Bbb{C%
}\left[ \left[ \zeta \right] \right] }{\frak{A}}\rightarrow A^{p,q+1}\left(
X_{0}\right) \otimes \frac{\Bbb{C}\left[ \left[ \zeta \right] \right] }{%
\frak{A}}
\end{equation*}
is therefore considered as a deformation of the operator 
\begin{equation*}
\overline{\partial }_{\beta _{0}}:A^{p,q}\left( X_{0}\right) \rightarrow
A^{p,q+1}\left( X_{0}\right)
\end{equation*}
where 
\begin{equation*}
L_{0}=L_{\beta _{0}}.
\end{equation*}
So we always compute on the same space of \textit{ordinary} differential
forms on $X_{0}$ but simply deform the \textit{twisted} operator $\overline{%
\partial }_{\beta _{0}}$ via its formal deformation $\overline{\partial }%
_{\beta }$.

For any form $\varepsilon \in A^{p,q}\left( X_{0}\right) \otimes \Bbb{C}%
\left[ \left[ \zeta \right] \right] $ we denote the $m$-th summand under
this decomposition as 
\begin{equation*}
\varepsilon _{m}.
\end{equation*}
In the following we will use the identification $\left( \ref{6.0}\right) $
throughout. Now consider the $\Bbb{R}$-linear inclusion 
\begin{equation*}
H^{1}\left( \mathcal{O}_{X_{0}}\right) \rightarrow H^{1}\left( X_{0};\Bbb{R}%
\right) \subseteq H^{1}\left( X_{0};\Bbb{C}\right)
\end{equation*}
whose composition with $\left( \ref{3.10'}\right) $ is the identity. Then
the elements $\beta \in \mathrm{Hom}_{\Bbb{R}}\left( H_{1}\left( X_{0};\Bbb{Z%
}\right) ,\Bbb{R}\right) $ constructed above depend linearly on the
parameter $\zeta $ so that, abusing notation by denoting the linear family $%
\beta \left( \zeta \right) $ of linear functionals simply as $\beta $ we
have that 
\begin{equation*}
\overline{\partial }\beta \in A^{0,1}\left( X_{0}\right) \otimes \Bbb{C}%
\left[ \left[ \zeta \right] \right]
\end{equation*}
has pure weight one so that 
\begin{equation*}
\overline{\partial }\beta =\left( \overline{\partial }\beta \right) _{1}.
\end{equation*}

We will fix a line bundle 
\begin{equation*}
L_{0}=L_{\beta _{0}}
\end{equation*}
for some fixed 
\begin{equation*}
\beta _{0}\in H^{1}\left( X_{0};\Bbb{C}\right)
\end{equation*}
and let 
\begin{eqnarray*}
\partial _{0} &=&\partial -\left( \partial \beta _{0}\wedge \right) \\
\overline{\partial }_{0} &=&\overline{\partial }-\left( \overline{\partial }%
\beta _{0}\wedge \right) .
\end{eqnarray*}
Then writing 
\begin{equation*}
H^{*}\left( \Omega _{X_{0}}^{p},\overline{\partial }_{0}\right) :=H\left(
A_{X_{0}}^{p,*},\overline{\partial }_{0}\right)
\end{equation*}
we have that 
\begin{equation*}
\mathcal{H}^{q}\left( H^{*}\left( \Omega _{X_{0}}^{p},\overline{\partial }%
_{0}\right) \otimes \frac{\Bbb{C}\left[ \left[ \zeta \right] \right] }{\frak{%
A}},\overline{\partial }\beta \wedge \right)
\end{equation*}
has a natural structure as a graded module over 
\begin{equation*}
\Delta _{\frak{A}}=Spec\frac{\Bbb{C}\left[ \left[ \zeta \right] \right] }{%
\frak{A}}.
\end{equation*}

\begin{lemma}
\label{l6}Let $\frak{A}$ be a homogeneous ideal in $\Bbb{C}\left[ \left[
\zeta \right] \right] $. Then over $\Delta _{\frak{A}}$ we have the
following:

i) For a solution 
\begin{equation*}
\overline{\partial }_{0}\left( \alpha _{n}+\alpha _{n+1}\right) \equiv 
\overline{\partial }\beta \wedge \left( \alpha _{n}+\alpha _{n+1}\right)
\end{equation*}
modulo $\frak{m}^{n+2}+\frak{A}$, 
\begin{equation}
\begin{array}{c}
\overline{\partial }_{0}\alpha _{n}=0 \\ 
\left\{ \overline{\partial }\beta \wedge \alpha _{n}\right\} =0\in
H^{p,q+1}\left( X_{0}\right) \otimes \frac{\frak{m}^{n+1}+\frak{A}}{\frak{m}%
^{n+2}+\frak{A}}
\end{array}
\label{6.10'}
\end{equation}
and each element $\alpha _{n}$ satisfying $\left( \ref{6.10'}\right) $
extends to a solution 
\begin{equation*}
\alpha _{n}+\alpha _{n+1}
\end{equation*}
modulo $\frak{m}^{n+2}+\frak{A}$, where 
\begin{equation*}
\overline{\partial }_{0}\alpha _{n+1}=\overline{\partial }\beta \wedge
\alpha _{n}.
\end{equation*}
ii) A solution 
\begin{equation*}
\alpha _{n}+\alpha _{n+1}
\end{equation*}
modulo $\frak{m}^{n+2}+\frak{A}$ has the form 
\begin{equation*}
\left( \overline{\partial }_{0}-\overline{\partial }\beta \wedge \right)
\left( \varepsilon _{n}+\varepsilon _{n+1}\right)
\end{equation*}
for some $\left( \varepsilon _{n}+\varepsilon _{n+1}\right) $ if and only if 
\begin{eqnarray*}
\overline{\partial }_{0}\varepsilon _{n}=\alpha _{n} \\
\alpha _{n+1}=\overline{\partial }\beta \wedge \varepsilon _{n}+\overline{%
\partial }_{0}\varepsilon _{n+1}.
\end{eqnarray*}

iii) For any solution 
\begin{equation*}
\overline{\partial }_{0}\left( \alpha _{n}+\alpha _{n+1}\right) \equiv 
\overline{\partial }\beta \wedge \left( \alpha _{n}+\alpha _{n+1}\right)
\end{equation*}
modulo $\frak{m}^{n+2}+\frak{A}$ with 
\begin{equation*}
\partial _{0}\alpha _{n}=0
\end{equation*}
there is a solution 
\begin{equation*}
\overline{\partial }_{0}\left( \alpha _{n}+\partial \gamma _{n+1}\right)
\equiv \overline{\partial }\beta \wedge \left( \alpha _{n}+\partial
_{0}\gamma _{n+1}\right) .
\end{equation*}
Also 
\begin{equation*}
\overline{\partial }_{0}\left( \alpha _{n+1}-\partial _{0}\gamma
_{n+1}\right) =0.
\end{equation*}

iv) For $n>0$, any $n$-order solution 
\begin{equation}
\overline{\partial }_{0}\left( \alpha _{0}+\ldots +\alpha _{n-1}+\partial
_{0}\gamma _{n}\right) \equiv \overline{\partial }\beta \wedge \left( \alpha
_{0}+\ldots +\alpha _{n-1}+\partial _{0}\gamma _{n}\right)  \notag
\end{equation}
modulo $\frak{m}^{n+1}+\frak{A}$ with $\alpha _{0}\in $ $A^{p,q}\left(
X_{0}\right) $ and $\gamma _{i}\in $ $A^{p-1,q}\left( X_{0}\right) $ extends
to an $\left( n+1\right) $-st order solution 
\begin{equation*}
\overline{\partial }_{0}\left( \alpha _{0}+\ldots +\alpha _{n-1}+\partial
_{0}\gamma _{n}+\partial _{0}\gamma _{n+1}\right) \equiv \overline{\partial }%
\beta \wedge \left( \alpha _{0}+\ldots +\alpha _{n-1}+\partial _{0}\gamma
_{n}+\partial _{0}\gamma _{n+1}\right)
\end{equation*}
modulo $\frak{m}^{n+2}+\frak{A}$.
\end{lemma}

\begin{proof}
i) and ii) are obvious.

iii) Since 
\begin{equation*}
\overline{\partial }_{0}\alpha _{n+1}=\overline{\partial }\beta \wedge
\alpha _{n}
\end{equation*}
we have that 
\begin{equation*}
\overline{\partial }\beta \wedge \alpha _{n}
\end{equation*}
is $\partial _{0}$-closed and $\overline{\partial }_{0}$-exact and so by the
Hodge $\partial \overline{\partial }$-lemma in Theorem \ref{p5}ii) there
exist $\gamma _{n+1}$ such that 
\begin{equation*}
\overline{\partial }\beta \wedge \alpha _{n}=\overline{\partial }
_{0}\partial _{0}\gamma _{n+1}.
\end{equation*}
Thus 
\begin{equation*}
\overline{\partial }_{0}\left( \alpha _{n}+\partial _{0}\gamma _{n+1}\right)
\equiv \overline{\partial }\beta \wedge \left( \alpha _{n}+\partial
_{0}\gamma _{n+1}\right) .
\end{equation*}
Clearly 
\begin{equation*}
\overline{\partial }_{0}\left( \alpha _{n+1}-\partial _{0}\gamma
_{n+1}\right) =0.
\end{equation*}

iv) The form 
\begin{equation*}
\overline{\partial }\beta \wedge \partial _{0}\gamma _{n}
\end{equation*}
is $\partial _{0}$-exact and $\overline{\partial }_{0}$-closed and so by the
Hodge $\partial \overline{\partial }$-lemma (Theorem \ref{p5}ii)) 
\begin{equation*}
\overline{\partial }\beta \wedge \partial _{0}\gamma _{n}=\overline{\partial 
}_{0}\partial _{0}\gamma _{n+1}
\end{equation*}
for some $\gamma _{n+1}\in $ $A^{p-1,q}\left( X_{0}\right) \otimes \frac{%
\frak{m}^{n+1}+\frak{A}}{\frak{m}^{n+2}+\frak{A}}.$
\end{proof}

We use this last lemma in the proof of the following:

\begin{theorem}
\label{c6.1}i) Let $\frak{A}$ be a homogeneous ideal centered at a point 
\begin{equation*}
\left\{ L_{0}\right\} \in \mathrm{Pic}^{0}\left( X_{0}\right) .
\end{equation*}
Then with respect to the filtration of 
\begin{equation*}
R^{q}\pi _{*}\left( \Omega _{X_{0}}^{p}\left( L_{0}\right) \otimes \left.
a^{*}P\right| _{\Delta _{\frak{A}}}\right)
\end{equation*}
by the subspaces generated by 
\begin{equation*}
R^{q}\pi _{*}\left( \Omega _{X_{0}}^{p}\left( L_{0}\right) \otimes \left.
a^{*}\frak{m}^{n}P\right| _{\Delta _{\frak{A}}}\right) ,
\end{equation*}
we have 
\begin{eqnarray}
Gr\left( R^{q}\pi _{*}\left( \Omega _{X_{0}}^{p}\left( L_{0}\right) \otimes
\left. a^{*}P\right| _{\Delta _{\frak{A}}}\right) \right) &\cong &\mathcal{H}%
^{q}\left( H^{*}\left( \Omega _{X_{0}}^{p},\overline{\partial }_{0}\right)
\otimes \frac{\Bbb{C}\left[ \left[ \zeta \right] \right] }{\frak{A}},\ 
\overline{\partial }\beta \wedge \right)  \label{new6.2} \\
&\cong &\mathcal{H}^{q}\left( H^{*}\left( \Omega _{X_{0}}^{p}\left(
L_{0}\right) \right) \otimes \frac{\Bbb{C}\left[ \left[ \zeta \right]
\right] }{\frak{A}},\ \overline{\partial }\beta \wedge \right) .
\end{eqnarray}
ii) There exists a (non-canonical) filtration preserving isomorphism of $%
\mathcal{O}_{\Delta _{\frak{A}}}$-modules 
\begin{equation}
R^{q}\pi _{*}\left( \Omega _{X_{0}}^{p}\left( L_{0}\right) \otimes \left.
a^{*}P\right| _{\Delta _{\frak{A}}}\right) \rightarrow Gr\left( R^{q}\pi
_{*}\left( \Omega _{X_{0}}^{p}\left( L_{0}\right) \otimes \left.
a^{*}P\right| _{\Delta _{\frak{A}}}\right) \right) .  \label{6.3}
\end{equation}

(Compare \cite{GL2}, Theorem 3.2.)
\end{theorem}

\begin{proof}
i) By $\left( \ref{3.3}\right) $ 
\begin{equation*}
R^{q}\pi _{*}\left( \Omega _{X_{0}}^{p}\otimes \left. a^{*}P\right| _{\Delta
_{\frak{A}}}\right) \cong \mathcal{H}^{q}\left( A_{X_{0}}^{p,*}\otimes \frac{%
\Bbb{C}\left[ \left[ \zeta \right] \right] }{\frak{A}},\overline{\partial }
_{0}-\overline{\partial }\beta \wedge \right) .
\end{equation*}
We filter 
\begin{equation*}
A_{X_{0}}^{p,*}\otimes \frac{\Bbb{C}\left[ \left[ \zeta \right] \right] }{%
\frak{A}}
\end{equation*}
by subcomplexes 
\begin{equation*}
C_{n}:=A_{X_{0}}^{p,*}\otimes \frac{\frak{m}^{n}+\frak{A}}{\frak{A}}
\end{equation*}
and form the associated spectral sequence with 
\begin{equation*}
E_{1}^{n+q,-n}=H^{q}\left( \Omega _{X_{0}}^{p},\overline{\partial }
_{0}\right) \otimes \frac{\frak{m}^{n}+\frak{A}}{\frak{m}^{n+1}+\frak{A}}
\end{equation*}
and 
\begin{equation*}
d_{1}:E_{1}^{n+q,-n}\rightarrow E_{1}^{n+q+2,-n-1}
\end{equation*}
given by 
\begin{equation*}
d_{1}\left( \alpha \right) =-\overline{\partial }_{0}\beta \wedge \alpha .
\end{equation*}
By Lemma \ref{l6}i-ii), assertion i) of the theorem will follow from the
fact that this spectral sequence degenerates at $E_{2}$, that is, 
\begin{equation*}
d_{r}=0
\end{equation*}
for all $r\geq 2$. To compute $d_{2}$ let 
\begin{equation*}
Z_{n}^{p,q}=\ker \left( \overline{\partial }_{0}:A^{p,q}\left( X_{0}\right)
\otimes \frac{\frak{m}^{n}+\frak{A}}{\frak{m}^{n+1}+\frak{A}}\rightarrow
A^{p,q+1}\left( X_{0}\right) \otimes \frac{\frak{m}^{n}+\frak{A}}{\frak{m}
^{n+1}+\frak{A}}\right) .
\end{equation*}
and let $K_{n}^{p,q}$ denote the intersection 
\begin{equation*}
\ker \left( \overline{\partial }\beta \wedge :Z_{n}^{p,q}\rightarrow
H^{q+1}\left( \Omega _{X_{0}}^{p},\overline{\partial }_{0}\right) \otimes 
\frac{\frak{m}^{n+1}+\frak{A}}{\frak{m}^{n+2}+\frak{A}}\right) \cap \ker
\left( \left. \partial _{0}\right| _{Z_{n}^{p,q}}\right) .
\end{equation*}

By Lemma \ref{l6}iii-iv) we can define a $\Bbb{C}$-linear function 
\begin{eqnarray}
K_{n}^{p,q} &\rightarrow &\left( \left. \mathcal{H}^{q}\left(
A_{X_{0}}^{p,*}\otimes \frac{\frak{m}^{n}+\frak{A}}{\frak{A}},\overline{%
\partial }_{0}-\overline{\partial }\beta \wedge \right) \right| _{\Delta _{%
\frak{A}}}\right) .  \label{6.01} \\
\alpha _{n} &\mapsto &\left( \alpha _{n}+\partial _{0}\gamma _{n+1}+\ldots
\right)  \notag
\end{eqnarray}
and the natural map 
\begin{equation*}
K_{n}^{p,q}\rightarrow E_{2}^{n+q,-n}
\end{equation*}
is surjective. So we can compute $d_{2}$ on an element $\eta \in
E_{2}^{n-q,-n}$ by applying the operator $\overline{\partial }_{0}-\overline{%
\partial }\beta \wedge $ to a preimage $\tilde{\eta}\in K_{n}^{p,q}$. So we
can compute $d_{2}$ on an element $\eta \in E_{2}^{n-q,-n}$ by applying the
operator $\overline{\partial }_{0}-\overline{\partial }\beta \wedge $ to a
preimage $\tilde{\eta}\in K_{n}^{p,q}$. But using Lemma \ref{l6}iv)
repeatedly we can recursively pick our representative $\tilde{\eta}$ for the
equivalence class $\eta $ such that 
\begin{equation*}
\left( \overline{\partial }_{0}-\overline{\partial }\beta \wedge \right)
\left( \tilde{\eta}\right) =0
\end{equation*}
in $A^{p,q+1}\left( X_{0}\right) \otimes \frac{\frak{m}^{n+2}+\frak{A}}{%
\frak{m}^{n+2+r}+\frak{A}}$ for all values of $r\geq 2$. Thus 
\begin{equation*}
d_{r}=0
\end{equation*}
for all $r\geq 2$.

ii) We wish to construct a mapping 
\begin{equation*}
E_{\infty }^{n+q,-n}\rightarrow \mathcal{H}^{q}\left( A_{X_{0}}^{p,*}\otimes 
\frac{\frak{m}^{n}+\frak{A}}{\frak{A}},\overline{\partial }_{0}-\overline{
\partial }\beta \wedge \right)
\end{equation*}
whose composition with the natural map 
\begin{equation*}
\mathcal{H}^{q}\left( A_{X_{0}}^{p,*}\otimes \frac{\frak{m}^{n}+\frak{A}}{%
\frak{A}},\overline{\partial }_{0}-\overline{\partial }\beta \wedge \right)
\rightarrow E_{\infty }^{n+q,-n}
\end{equation*}
is the identity. The mapping $\left( \ref{6.01}\right) $ above can be chosen
in such a way as to do the job. Namely, if $\alpha _{n}$ is in the kernel of
the composition of $\left( \ref{6.01}\right) $ and the projection 
\begin{equation*}
\left( \left. \mathcal{H}^{q}\left( A_{X_{0}}^{p,*}\otimes \frac{\frak{m}
^{n}+\frak{A}}{\frak{A}},\overline{\partial }_{0}-\overline{\partial }\beta
\wedge \right) \right| _{\Delta _{\frak{A}}}\right) \rightarrow
E_{2}^{n+q,-n}
\end{equation*}
then 
\begin{equation*}
\alpha _{n}=\overline{\partial }_{0}\varepsilon _{n}-\overline{\partial }
\beta \wedge \varepsilon _{n-1}
\end{equation*}
where 
\begin{equation*}
\overline{\partial }_{0}\varepsilon _{n-1}=0.
\end{equation*}
So we can map $\alpha _{n}$ to 
\begin{equation*}
\left( \overline{\partial }_{0}-\overline{\partial }\beta \wedge \right)
\left( \varepsilon _{n-1}+\varepsilon _{n}\right)
\end{equation*}
which lies in 
\begin{equation*}
\mathcal{H}^{q}\left( A_{X_{0}}^{p,*}\otimes \frac{\frak{m}^{n}+\frak{A}}{%
\frak{A}},\overline{\partial }_{0}-\overline{\partial }\beta \wedge \right)
\end{equation*}
and is zero in 
\begin{equation*}
\mathcal{H}^{q}\left( A_{X_{0}}^{p,*}\otimes \frac{\Bbb{C}\left[ \left[
\zeta \right] \right] }{\frak{A}},\overline{\partial }_{0}-\overline{
\partial }\beta \wedge \right) .
\end{equation*}
\end{proof}

\begin{corollary}
\label{c6.2}The schemes 
\begin{equation*}
h^{q}\left( \Omega _{X_{0}}^{p}\otimes L\right) \geq m_{p,q}
\end{equation*}
are finite unions of (reduced) linear subspaces at $\left\{ L_{0}\right\}
\in \mathrm{Pic}^{0}\left( X_{0}\right) $ with tangent cones computed by 
\begin{equation*}
\frac{\ker \left( H^{q}\left( \Omega _{X_{0}}^{p}\left( L_{0}\right) ,%
\overline{\partial }\right) \otimes \frac{\frak{m}}{\frak{m}^{2}}\overset{%
\overline{\partial }\beta \wedge }{\longrightarrow }H^{q+1}\left( \Omega
_{X_{0}}^{p}\left( L_{0}\right) ,\overline{\partial }\right) \otimes \frac{%
\frak{m}^{2}}{\frak{m}^{3}}\right) }{\mathrm{im}\left( H^{q-1}\left( \Omega
_{X_{0}}^{p}\left( L_{0}\right) ,\overline{\partial }\right) \overset{%
\overline{\partial }\beta \wedge }{\longrightarrow }H^{q}\left( \Omega
_{X_{0}}^{p}\left( L_{0}\right) ,\overline{\partial }\right) \otimes \frac{%
\frak{m}}{\frak{m}^{2}}\right) }.
\end{equation*}
\end{corollary}

\begin{proof}
Restrict $\beta \left( \zeta \right) $ to a one-dimensional linear subspace $%
S$ of $H^{1}\left( \mathcal{O}_{X_{0}}\right) .$ If 
\begin{equation*}
h^{q}\left( \Omega _{X_{0}}^{p}\otimes L_{0}\otimes L_{\beta }\right) \geq
m_{p,q}
\end{equation*}
to first-order at $0\in S$, then the same is true to all orders on an
analytic neighborhood of $0\in S$.
\end{proof}

Finally we give a relative version of the above results. Let 
\begin{equation*}
p^{\prime }:X\rightarrow X^{\prime }
\end{equation*}
be a smooth family of compact K\"{a}hler manifolds over a smooth parameter
space $X^{\prime }$ and let 
\begin{equation*}
\left\{ L_{0}^{X^{\prime }}\right\}
\end{equation*}
be the image of a section of 
\begin{equation*}
\sigma :X^{\prime }\rightarrow \mathrm{Pic}^{0}\left( X/X^{\prime }\right) .
\end{equation*}
Then under the automorphism of $\mathrm{Pic}^{0}\left( X/X^{\prime }\right) $
induced by translation by $\sigma $ we have an isomorphism between 
\begin{equation*}
\frak{S}:=\widetilde{Sym}\left( \left( R^{1}p_{*}^{\prime }\mathcal{O}%
_{X}\right) ^{\vee }\right)
\end{equation*}
(where $\widetilde{Sym}$ denotes the completion of the symmetric algebra at
the zero section) and the completion of ring of regular functions on $%
\mathrm{Pic}^{0}\left( X/X^{\prime }\right) $ with respect to the ideal of $%
\sigma \left( X^{\prime }\right) .$ Then for 
\begin{equation*}
\pi :X\times _{X^{\prime }}\mathrm{Pic}^{0}\left( X/X^{\prime }\right)
\rightarrow \mathrm{Pic}^{0}\left( X/X^{\prime }\right) ,
\end{equation*}
the relative Poincar\'{e} bundle 
\begin{equation*}
P_{X/X^{\prime }}\rightarrow \mathrm{Alb}\left( X/X^{\prime }\right) \times
_{X^{\prime }}\mathrm{Pic}^{0}\left( X/X^{\prime }\right) ,
\end{equation*}
and isomorphism 
\begin{equation*}
\beta :R^{1}p_{*}^{\prime }\left( \mathcal{O}_{X/X^{\prime }}\right)
\rightarrow \mathrm{Hom}_{\Bbb{R}}\left( p_{*}^{\prime }\left( \Omega
_{X/X^{\prime }}^{1}\right) ^{\vee },\Bbb{R}\right)
\end{equation*}
constructed in \S \ref{3} and \S \ref{6} we have:

\begin{corollary}
\begin{equation*}
R^{q}\pi _{*}\left( \Omega _{X/X^{\prime }}^{p}\left( L_{0}\right) \otimes
\left. a^{*}P_{X/X^{\prime }}\right| _{\Delta _{\frak{A}}}\right) \cong 
\mathcal{H}^{q}\left( R^{*}p_{*}^{\prime }\left( \Omega _{X/X^{\prime
}}^{p}\left( L_{0}\right) \right) \otimes _{\mathcal{O}_{X^{\prime }}}\frac{%
\mathcal{O}_{X^{\prime }}\left[ \left[ \zeta \right] \right] }{\frak{A}},%
\overline{\partial }\beta \wedge \right) .
\end{equation*}
\end{corollary}

\begin{proof}
Harmonic decomposition of relative forms over $X^{\prime }$ and all parts of
Lemma \ref{l6} and of Theorem \ref{c6.1} can be carried out over $X^{\prime
} $.
\end{proof}

\section{The twisted coefficients case\label{7}}

Let 
\begin{equation*}
\sum B_{j}
\end{equation*}
be a simple-normal-crossing divisor on $X_{0}$, 
\begin{equation*}
B=\sum b_{j}B_{j}
\end{equation*}
for $b_{j}$ positive integers, and let $M$ be a line bundle on $X_{0}$ such
that 
\begin{equation*}
M^{\otimes N}=\mathcal{O}_{X_{0}}\left( B\right)
\end{equation*}
for some integers $N,\ b_{j}>0$. Following \S 3 of \cite{EV1} we will write 
\begin{equation}
{M^{(i)}}:=M^{\otimes i}\otimes \mathcal{O}_{X_{0}}(-\sum \lfloor \frac{
ib_{j}}{N}\rfloor B_{j}),  \label{7.1}
\end{equation}
and 
\begin{equation}
{B^{(i)}}:=\sum_{N\ \nmid \ ib_{j}}B_{j}.  \label{7.2}
\end{equation}
Let 
\begin{equation*}
f:Y_{0}\rightarrow X_{0}
\end{equation*}
be the normalization of the $N$-cyclic cover given by the inverse image of
the canonical section $s$ of $\mathcal{O}_{X_{0}}\left( \sum
b_{j}B_{j}\right) $ under the $N$-th power map 
\begin{equation*}
{M\rightarrow }M^{\otimes N}.
\end{equation*}
Then $Y_{0}$ is nonsingular. In fact in \S 3 of \cite{EV1} it is shown that
the natural algebra structure on 
\begin{equation*}
\frac{\bigoplus_{i=0}^{\infty }({M^{(i)}})^{\vee }}{\left\{ \alpha \in {\
M^{-N}}\sim s\left( \alpha \right) \in \mathcal{O}_{X_{0}}\right\} }
\end{equation*}
has spectrum $Y_{0}$ so that 
\begin{equation}
f_{*}\mathcal{O}_{Y_{0}}=\bigoplus_{i=0}^{N-1}({M^{(i)}})^{\vee },
\label{7.3}
\end{equation}
and furthermore 
\begin{equation}
f_{*}\Omega _{Y_{0}}^{p}=\bigoplus_{i=0}^{N-1}\Omega _{X_{0}}^{p}(log{B^{(i)}%
})\otimes ({M^{(i)}})^{\vee }.  \label{7.4}
\end{equation}
We remark that the above covering is Galois with cyclic Galois group of
order $N$. The above decompositions are simply the decomposition into direct
sums of sheaves of eigenvectors. For $i=0$ one has $\Omega _{X_{0}}^{p}(log{%
\ \ \ B^{(0)}})\otimes ({M^{(0)}})^{\vee }=\Omega _{X_{0}}^{p}$.

\begin{theorem}
\label{c7.1}Let $\frak{A}$ be a homogeneous ideal as in \S 6 and let $%
\left\{ L_{0}\right\} \in \mathrm{Pic}^{0}\left( X_{0}\right) $. Then with
respect to the filtration of 
\begin{equation*}
R^{q}\pi _{*}\left( \Omega _{X_{0}}^{p}(log{B^{(i)}})\otimes ({M^{(i)}}%
)^{\vee }\otimes L_{0}\otimes \left. a^{*}P\right| _{\Delta _{\frak{A}%
}}\right)
\end{equation*}
by the subspaces generated by 
\begin{equation*}
R^{q}\pi _{*}\left( \Omega _{X_{0}}^{p}(log{B^{(i)}})\otimes ({M^{(i)}}%
)^{\vee }\otimes L_{0}\otimes \left. a^{*}\frak{m}^{n}P\right| _{\Delta _{%
\frak{A}}}\right) ,
\end{equation*}
for $0\leq i\leq N-1$, we have 
\begin{eqnarray*}
&&Gr\left( R^{q}\pi _{*}\left( \Omega _{X_{0}}^{p}(log{B^{(i)}})\otimes ({\
M^{(i)}})^{\vee }\otimes L_{0}\otimes \left. a^{*}P\right| _{\Delta _{\frak{A%
}}}\right) \right) \\
&\cong &\mathcal{H}^{q}\left( H^{*}\left( \Omega _{X_{0}}^{p}(log{\ B^{(i)}}%
)\otimes L_{0}\otimes ({M^{(i)}})^{\vee }\right) \otimes \frac{\Bbb{C}\left[
\left[ \zeta \right] \right] }{\frak{A}},\ \overline{\partial }\beta \wedge
\right) .
\end{eqnarray*}
\end{theorem}

\begin{proof}
Let $\tilde{\pi }=\pi \circ (f\times id_{\Delta }):Y_{0}\times \Delta
\longrightarrow \Delta $. By Theorem \ref{c6.1} we have 
\begin{equation*}
Gr\left( R^{q}\tilde{\pi }_{*}\left( \Omega _{Y_{0}}^{p}\otimes
f^{*}L_{0}\otimes \left. (p\circ a)^{*}P\right| _{\Delta _{\frak{A}}}\right)
\right) \cong \mathcal{H}^{q}\left( H^{*}\left( \Omega _{Y_{0}}^{p},%
\overline{\partial }_{0}\right) \otimes \frac{\Bbb{C}\left[ \left[ \zeta
\right] \right] }{\frak{A}},\ f^{*}\overline{\partial }\beta \wedge \right) .
\end{equation*}
But the left-hand expression is isomorphic to 
\begin{equation*}
Gr\left( R^{q}{\pi }_{*}\left( \bigoplus_{i=0}^{N-1}\Omega _{X_{0}}^{p}(log{%
\ \ B^{(i)}})\otimes ({M^{(i)}})^{\vee }\otimes L_{0}\otimes \left.
a^{*}P\right| _{\Delta _{\frak{A}}}\right) \right) .
\end{equation*}
And the right-hand expression is isomorphic to 
\begin{equation*}
\mathcal{H}^{q}\left( \bigoplus_{i=0}^{N-1}H^{*}\left( \Omega
_{X_{0}}^{p}(log{B^{(i)}})\otimes L_{0}\otimes ({M^{(i)}})^{\vee }\right)
\otimes \frac{\Bbb{C}\left[ \left[ \zeta \right] \right] }{\frak{A}},\ 
\overline{\partial }\beta \wedge \right) .
\end{equation*}
Comparing eigenspaces, the theorem now follows.
\end{proof}

\begin{corollary}
\label{c7.2}The schemes 
\begin{equation}
h^{q}\left( \Omega _{X_{0}}^{p}(log{B^{(i)}})\otimes L\otimes ({M^{(i)}}%
)^{\vee }\right) \geq m_{p,q}  \label{lb''}
\end{equation}
are finite unions of (reduced) linear subspaces at $\left\{ L_{0}\right\}
\in \mathrm{Pic}^{0}\left( X_{0}\right) $ with tangent cones computed by the
cohomology of 
\begin{equation*}
\begin{array}{c}
H^{q-1}\left( \Omega _{X_{0}}^{p}(log{B^{(i)}})\otimes L_{0}\otimes ({M^{(i)}%
})^{\vee },\overline{\partial }\right) \\ 
\downarrow \overline{\partial }\beta \wedge \\ 
H^{q}\left( \Omega _{X_{0}}^{p}(log{B^{(i)}})\otimes L_{0}\otimes ({M^{(i)}}%
)^{\vee },\overline{\partial }\right) \otimes \frac{\frak{m}}{\frak{m}^{2}}
\\ 
\downarrow \overline{\partial }\beta \wedge \\ 
H^{q+1}\left( \Omega _{X_{0}}^{p}(log{B^{(i)}})\otimes L_{0}\otimes ({M^{(i)}%
})^{\vee },\overline{\partial }\right) \otimes \frac{\frak{m}^{2}}{\frak{m}%
^{3}}
\end{array}
.
\end{equation*}
In particular, in some analytic neighborhood of $\left\{ L_{0}\right\} \in 
\mathrm{Pic}^{0}\left( X_{0}\right) ,$ the schemes $\left( \ref{lb''}\right) 
$ are (reduced) unions of abelian subvarieties.
\end{corollary}

We remark that the results of \S 5 extend without difficulty to the case of
twisted coefficients as given in \S 7. Namely from \S 5 suppose that 
\begin{equation*}
\beta \in H^{1}\left( X_{0};\Bbb{R}\right) .
\end{equation*}
Then for 
\begin{equation*}
L=L_{\beta }
\end{equation*}
we have that 
\begin{equation*}
\overline{L}=L^{\vee }
\end{equation*}
and for flat line bundles $\tilde{L}=f^{*}L$ on $Y_{0}$ we have 
\begin{equation*}
\overline{H^{q}\left( \Omega _{Y_{0}}^{p}\otimes \tilde{L}\right) }\cong
H^{p}\left( \Omega _{Y_{0}}^{q}\otimes \tilde{L}^{\vee }\right)
\end{equation*}
and 
\begin{equation*}
\sum\nolimits_{p+q=r}H^{q}\left( \Omega _{Y_{0}}^{p}\otimes \tilde{L}\right)
\oplus \overline{H^{q}\left( \Omega _{Y_{0}}^{p}\otimes \tilde{L}\right) }%
\cong H^{r}\left( Y_{0};\tilde{L}_{\Bbb{R}}\right) \otimes \Bbb{C}
\end{equation*}
But recalling that 
\begin{equation*}
{M^{(i)}}=M^{\otimes i}\otimes \mathcal{O}_{X_{0}}(-\lfloor \frac{i}{N}%
B\rfloor )
\end{equation*}
and comparing eigenspace decompositions 
\begin{eqnarray*}
H^{q}\left( \Omega _{Y_{0}}^{p}\otimes \tilde{L}\right) &=&H^{q}\left(
f_{*}\Omega _{Y_{0}}^{p}\otimes L\right) \\
&=&\bigoplus_{i=0}^{N-1}H^{q}\left( \Omega _{X_{0}}^{p}(log{B^{(i)}})\otimes
({M^{(i)}})^{\vee }\otimes L\right)
\end{eqnarray*}
and 
\begin{eqnarray*}
H^{p}\left( \Omega _{Y_{0}}^{q}\otimes \tilde{L}^{\vee }\right)
&=&H^{p}\left( f_{*}\Omega _{Y_{0}}^{q}\otimes L^{\vee }\right) \\
&=&\bigoplus_{i=0}^{N-1}H^{p}\left( \Omega _{X_{0}}^{q}(log{B^{(i)}})\otimes
({M^{(i)}})^{\vee }\otimes L^{\vee }\right)
\end{eqnarray*}
we see that the $\Bbb{C}$-antilinear isomorphisms which generalize Hodge
symmetry in Theorem \ref{p5} are 
\begin{equation}
H^{q}\left( \Omega _{X_{0}}^{p}(log{B^{(i)}})\otimes \left( {M^{(i)}}\right)
^{\vee }\otimes L\right) \longrightarrow H^{p}\left( \Omega _{X_{0}}^{q}(log{%
\ B^{(N-i)}})\otimes \left( {M^{(N-i)}}\right) ^{\vee }\otimes L^{\vee
}\right) .  \label{7.5}
\end{equation}
Note that $\log B^{(N-i)}=\log B^{(i)}$ and 
\begin{eqnarray}
{M^{(i)}\otimes }M^{(N-i)} &=&\mathcal{O}_{X_{0}}\left( B\right) \otimes 
\mathcal{O}_{X_{0}}(-\lfloor \frac{i}{N}B\rfloor )\otimes \mathcal{O}%
_{X_{0}}(-\lfloor \frac{N-i}{N}B\rfloor )  \label{7.6} \\
&=&\mathcal{O}_{X_{0}}\left( B^{(i)}\right) .  \notag
\end{eqnarray}
Thus the right hand side of $\left( \ref{7.5}\right) $ is isomorphic to 
\begin{equation*}
H^{p}(\Omega _{X_{0}}^{q}\left( logB^{(i)}\right) \otimes
M^{(i)}(-B^{(i)})\otimes L^{\vee })
\end{equation*}
as in the isomorphism of K. Timmerscheidt (see \cite{EV1}, Theorem 13.5).

In the following we will be interested in the line-bundle components 
\begin{equation}
f_{*}\omega _{Y_{0}}=\bigoplus_{i=0}^{N-1}\omega _{X_{0}}(log{B^{(i)}}
)\otimes ({M^{(i)}})^{\vee }  \label{7pfd}
\end{equation}
of the push-foward of the dualizing sheaf of $Y_{0}$. Since 
\begin{equation*}
\omega _{X_{0}}(log{B^{(i)}})=\omega _{X_{0}}({B^{(i)}}),
\end{equation*}
by $\left( \ref{7.6}\right) $ we can rewrite $\left( \ref{7pfd}\right) $ as 
\begin{equation}
f_{*}\omega _{Y_{0}}=\omega _{X_{0}}\oplus \left(
\bigoplus\nolimits_{i=1}^{N-1}M^{\left( i\right) }\right) .  \label{7pfd'}
\end{equation}

\section{Generic Vanishing Theorems}

From now on we will drop the subscripts and simply write $X$ in place of $%
X_{0}$ and $L$ in place of $L_{\beta }$, and we always assume that 
\begin{equation*}
\beta \in H^{1}\left( X;\Bbb{R}\right)
\end{equation*}
so that 
\begin{equation*}
\overline{L}=L^{\vee }.
\end{equation*}
Also throughout the remainder of the paper 
\begin{equation*}
n:=\dim X.
\end{equation*}
For any coherent sheaf $\frak{F}$ on $X$ we will write 
\begin{equation*}
V^{q}(\Omega _{X}^{p}\otimes \frak{F}):=\{L\in \mathrm{Pic}
^{0}(X):h^{q}(\Omega _{X}^{p}\otimes \frak{F}\otimes L)>0\}
\end{equation*}
We wish to illustrate a number of results that follow from Theorems \ref
{c6.1} and \ref{c7.1} concerning the geometry of certain loci $V^{q}(\Omega
_{X}^{p}\otimes \frak{F})$.

\begin{theorem}
\label{c8.2} Let $\left\{ L\right\} \in S$ be any general point of an
irreducible component $S$ of $V^{q}(\Omega _{X}^{p})$. Suppose that 
\begin{equation*}
v\in H^{1}(X,\mathcal{O}_{X})\cong T_{\left\{ L\right\} }\left( \mathrm{Pic}%
^{0}(X)\right)
\end{equation*}
is not tangent to $S$. Then the sequence 
\begin{equation*}
H^{q-1}(\Omega _{X}^{p}\left( L\right) )\overset{\cup v}{\longrightarrow }%
H^{q}(\Omega _{X}^{p}\left( L\right) )\overset{\cup v}{\longrightarrow }%
H^{q+1}(\Omega _{X}^{p}\left( L\right) )
\end{equation*}
is exact. If $\left\{ L\right\} $ is a general point of $S$ and $v$ is
tangent to $S$ at $\left\{ L\right\} $, then the maps in the above sequence
vanish.
\end{theorem}

\begin{proof}
When $p=0$ this is \cite{EL}, Theorem 1.2.3. The proof of the general case
follows immediately from Corollary \ref{c7.2}.
\end{proof}

Thus following Lemma 1.8 of \cite{EL} and \cite{GL1}, \cite{GL2} we have:

\begin{theorem}
\label{c8.3}

i) The loci $V^{0}(\Omega _{X}^{p})$, and so also the loci $V^{p}(\mathcal{O}%
_{X})$, are of dimension $\leq \left( p+\dim \left( \mathrm{Alb}\left(
X\right) \right) -\dim \left( a\left( X\right) \right) \right) .$

ii) Let $m$ denote $\dim \left( a\left( X\right) \right) $. Then 
\begin{equation*}
{\text{\textrm{Pic}}}^{0}\left( X\right) \supseteq V^{0}\left( X,\Omega
_{X}^{m}\right) \supseteq V^{0}\left( X,\Omega _{X}^{m-1}\right) \supseteq
...\supseteq V^{0}\left( X,\mathcal{O}_{X}\right) =\left\{ \mathcal{O}%
_{X}\right\} .
\end{equation*}
\end{theorem}

\begin{proof}
i) Let $\left\{ L\right\} $ be a general point of $V^{p}(\mathcal{O}_{X})$.
Then by Theorem \ref{c8.2} we have, for 
\begin{equation*}
W=\left\{ \partial \beta \in H^{0}\left( \Omega _{X}^{1}\right) :\overline{
\alpha }\wedge \overline{\partial \beta }=0\right\} ,
\end{equation*}
the inequality 
\begin{equation*}
\dim W\geq \dim \left( V^{p}(\mathcal{O}_{X})\right) .
\end{equation*}
Suppose 
\begin{equation*}
\alpha \left( L\right) \in A_{X}^{p,0}
\end{equation*}
is the harmonic representative for a non-zero element of 
\begin{equation*}
H^{0}\left( \Omega _{X}^{p}\left( L^{\vee }\right) \right) .
\end{equation*}
For general $x\in X$, let 
\begin{equation*}
W\left( x\right) =\left\{ \varepsilon \in \left. \Omega _{X}^{1}\right|
_{x}:\varepsilon \wedge \alpha \left( x\right) =0\right\} .
\end{equation*}
Now the codifferential of the Albanese map $a$ is given by 
\begin{equation*}
a_{x}^{*}:H^{0}\left( \Omega _{X}^{1}\right) \rightarrow \left. \Omega
_{X}^{1}\right| _{x}
\end{equation*}
so that, for a general point $x\in X$, $a_{x}^{*}$ has rank $=\dim \left(
a\left( X\right) \right) .$ So $\ker a_{x}^{*}$ has dimension 
\begin{equation*}
r:=\dim \left( \mathrm{Alb}X\right) -\dim \left( a\left( X\right) \right) .
\end{equation*}
Let 
\begin{equation*}
\gamma _{1},\ldots ,\gamma _{s}
\end{equation*}
be a basis of $W$, then we can assume that, for example, 
\begin{equation*}
\gamma _{1}\left( x\right) ,\ldots ,\gamma _{s-r}\left( x\right)
\end{equation*}
is a partial basis for the image of $W$ in $W\left( x\right) \subseteq
\left. \Omega _{X}^{1}\right| _{x}$. Writing the non-zero form $\alpha
\left( x\right) $ in terms of a completion of this partial basis to a full
basis of $\left. \Omega _{X}^{1}\right| _{x}$, the condition 
\begin{equation*}
\gamma _{i}\left( x\right) \wedge \alpha \left( x\right) =0
\end{equation*}
implies that 
\begin{equation*}
\alpha \left( x\right) =\gamma _{1}\left( x\right) \wedge \ldots \wedge
\gamma _{s-r}\left( x\right) \wedge \gamma ^{\prime }\left( x\right)
\end{equation*}
for some form $\gamma ^{\prime }\left( x\right) $. Thus $s-r\leq p$, that
is, 
\begin{equation*}
\dim W-\left( \dim \left( \mathrm{Alb}\left( X\right) \right) -\dim \left(
a\left( X\right) \right) \right) \leq p.
\end{equation*}
So 
\begin{equation*}
\dim \left( V^{p}(\mathcal{O}_{X})\right) \leq \dim W\leq p+\dim \left( 
\mathrm{Alb}\left( X\right) \right) -\dim \left( a\left( X\right) \right) .
\end{equation*}

ii) Suppose 
\begin{equation*}
H^{p}\left( L\right) \neq 0.
\end{equation*}
Then 
\begin{equation*}
H^{0}\left( \Omega _{X}^{p}\left( L^{\vee }\right) \right) \neq 0
\end{equation*}
so that there exists a non-trivial holomorphic $p$-form $\alpha $ with
coefficients in $L^{\vee }$. If $p<\dim a\left( X\right) $ there is a $%
\partial \beta $ such that 
\begin{equation*}
\partial \beta \wedge \alpha \neq 0
\end{equation*}
at the chain level and so also at the level of cohomology since holomorphic
forms are automatically harmonic. So 
\begin{equation*}
H^{0}\left( \Omega _{X}^{p+1}\left( L^{\vee }\right) \right) \neq 0
\end{equation*}
and therefore 
\begin{equation*}
H^{p+1}\left( L\right) \neq 0.
\end{equation*}
\end{proof}

By applying Theorem \ref{c8.3} to the case of a branched cover and using
Theorem \ref{c7.1} we have the following result that will be used in various
applications. Recall that in \S \ref{7} we let 
\begin{equation*}
\sum B_{j}
\end{equation*}
be a simple-normal-crossing divisor on $X$ and let $M$ be a line bundle on $%
X $ such that 
\begin{equation*}
M^{\otimes N}=\mathcal{O}_{X}\left( \sum b_{j}B_{j}\right)
\end{equation*}
for some integers $N,\ b_{j}>0$. We defined 
\begin{equation*}
{M^{(i)}}:=M^{\otimes i}\otimes \mathcal{O}_{X}(-\sum\nolimits_{j}\lfloor 
\frac{ib_{j}}{N}B_{j}\rfloor ).
\end{equation*}
Then from Theorem \ref{c7.1}, $\left( \ref{7pfd'}\right) $ and Theorem \ref
{c8.3} we conclude:

\begin{theorem}
\label{c8.4} i) Any irreducible component of 
\begin{equation*}
V^{q}({\omega }_{X}\otimes {M}^{\left( i\right) })
\end{equation*}
is a translate of a sub-torus of ${\text{\textrm{Pic}}}^{0}(X)$ and is of
codimension at least 
\begin{equation*}
q-({\text{\textrm{dim}}}(\mathrm{Alb}X)-{\text{\textrm{dim}}}(a\left(
X\right) )
\end{equation*}
in ${\text{\textrm{Pic}}}^{0}(X)$.

ii) Let $m$ denote $\dim \left( a\left( X\right) \right) $. Then for any $%
\left\{ L\right\} \in {\text{\textrm{Pic}}}^{0}(X)$%
\begin{equation*}
{\text{\textrm{Pic}}}^{0}(X)\supset V^{0}\left( {\omega }_{X}\otimes {M}%
^{\left( i\right) }\right) \supseteq V^{1}\left( {\omega }_{X}\otimes {M}%
^{\left( i\right) }\right) \supseteq \ldots \supseteq V^{m}\left( {\omega }%
_{X}\otimes {M}^{\left( i\right) }\right) .
\end{equation*}

iii) Let $\left\{ L\right\} \in S$ be a general point of an irreducible
component $S$ of $V^{q}(X,\Omega _{X}^{p}(log{B^{(i)}})\otimes ({M^{(i)}}%
)^{\vee })$. Suppose that 
\begin{equation*}
v\in H^{1}(X,\mathcal{O}_{X})\cong T_{\left\{ L\right\} }\left( {\text{%
\textrm{Pic}}}^{0}(X)\right)
\end{equation*}
is not tangent to $S$. Then the sequence 
\begin{equation*}
H^{q-1}(\Omega _{X}^{p}(log{B^{(i)}})\otimes ({M^{(i)}})^{\vee }\otimes L)%
\overset{\cup v}{\longrightarrow }H^{q}(\Omega _{X}^{p}(log{B^{(i)}})\otimes
({M^{(i)}})^{\vee }\otimes L)
\end{equation*}
\begin{equation*}
\overset{\cup v}{\longrightarrow }H^{q+1}(\Omega _{X}^{p}(log{B^{(i)}}%
)\otimes ({M^{(i)}})^{\vee }\otimes L)
\end{equation*}
is exact. If $\left\{ L\right\} $ is a general point of $S$, and $v$ is
tangent to $S$, then the maps in the above sequence vanish.
\end{theorem}

\section{Sections of line bundles on abelian varieties}

In this section we give an application of these methods to the zero sets of
sections of line bundles on an abelian variety. See \cite{EV2} for a related
statement on projective space. The methods are analogous to the ones of \cite
{E}, however the geometry of abelian varieties and Theorem \ref{c8.4} will
allow us to derive a somewhat more precise statement.

We will need the following proposition, analogous to \cite{H}, Proposition
2.1. Let 
\begin{equation*}
a:X\rightarrow a(X)\subset Alb (X)
\end{equation*}
be a generically finite morphism. 

Let $S$ be any irreducible component of $V^{0}\left( \omega _{X}\otimes
M^{(i)}\right) $. By Serre duality 
\begin{equation*}
h^{0}(\omega _{X}\otimes M^{(i)}\otimes L)=h^{n}((M^{(i)})^{\vee }\otimes
L^{\vee })
\end{equation*}
(where $n=dim(X)$). Therefore $S$ corresponds to a component of $%
V^{n}((M^{(i)})^{\vee })$ under the automorphism 
\begin{equation*}
\left\{ L\right\} \mapsto \left\{ L^{\vee }\right\}
\end{equation*}
of ${\text{\textrm{Pic}}}^{0}(X)$. $S$ is the translate of an abelian
subvariety of ${\text{\textrm{Pic}}}^{0}(X)$. The inclusion 
\begin{equation*}
S\rightarrow {\text{\textrm{Pic}}}^{0}(X)
\end{equation*}
induces a dual surjection of abelian varieties 
\begin{equation*}
d:Alb(X)\rightarrow S^{\vee }.
\end{equation*}
Suppose now that 
\begin{equation*}
B=\sum\nolimits_{j\geq 0}b_{j}B_{j}
\end{equation*}
has some component 
\begin{equation*}
B_{0}
\end{equation*}
such that the mapping 
\begin{equation*}
a:B_{0}\rightarrow a\left( B_{0}\right)
\end{equation*}
is generically finite, that is, equidimensional, and that 
\begin{equation*}
\left( d\circ a\right) \left( B_{0}\right) =S^{\vee }.
\end{equation*}
Suppose further that 
\begin{equation*}
N\nmid ib_{0}.
\end{equation*}

\begin{proposition}
\label{c9.1}For general $\left\{ L\right\} \in S\in {\text{\textrm{Pic}}}%
^{0}(X)$, all sections of the line bundle 
\begin{equation*}
\omega _{X}\otimes {M^{(i)}}\otimes a^{*}L
\end{equation*}
vanish on $B_{0}$.
\end{proposition}

\begin{proof}
Throughout this proof let 
\begin{equation*}
A:=\mathrm{Alb}\left( X\right)
\end{equation*}
Let $p$ be a general point of $B_{0}$. Let $U\subset A$ be an appropriate
neighborhood of $q:=a(p)$ such that 
\begin{equation*}
a(B_{0})\cap U
\end{equation*}
is smooth and defined on $\left( a(X)\cap U\right) $ by the equation $%
z_{1}=0 $. Let 
\begin{equation*}
w=\partial \alpha \in H^{0}(\Omega _{A}^{1})
\end{equation*}
be such that $\left( \partial \alpha \right) \left( q\right) =(d{z_{1}})_{q}$
. Since $d(a(B_{0}))=S^{\vee }$, it follows that 
\begin{equation*}
v=\overline{\partial \alpha }\notin T\left( S\right) \subset H^{1}(\mathcal{O%
}_{A}).
\end{equation*}
By Theorem \ref{c8.4}, the map 
\begin{equation*}
H^{n-1}\left( (M^{(i)})^{\vee }\otimes L^{\vee }\right) \overset{\cup v}{
\longrightarrow }H^{n}\left( (M^{(i)})^{\vee }\otimes L^{\vee }\right)
\end{equation*}
is surjective. So, by the $\Bbb{C}$-antilinear isomorphism (\ref{7.5}), the
map 
\begin{equation*}
H^{0}(\Omega _{X}^{n-1}(logB^{(i)})\otimes \left( M^{(i)}\left( -B^{\left(
i\right) }\right) \right) \otimes L)\overset{\wedge w}{\longrightarrow }
H^{0}(\omega _{X}\otimes M^{(i)}\otimes L)
\end{equation*}
is surjective. The assertion now follows from a local computation. Let $%
V=a^{-1}(U)$, we may assume that $B_{0}\cap V$ is also smooth and defined by
the equation $x_{1}=0$. Choosing local parameters $\{x_{1},...,x_{n}\}$ for $%
V$, we may write a section 
\begin{equation*}
s\in \Gamma (V;\Omega _{X}^{n-1}(logB^{(i)})\otimes M^{(i)}(-B^{(i)})\otimes
L)
\end{equation*}
locally in the form 
\begin{equation*}
s=\sum \eta _{k}x_{1}a_{k}
\end{equation*}
where 
\begin{equation*}
\eta _{1}=dx_{2}\wedge ...\wedge dx_{n},\ \ \ \eta _{k}=\frac{dx_{1}}{x_{1}}
\wedge ...\wedge dx_{k-1}\wedge dx_{k+1}\wedge ...\wedge dx_{n}
\end{equation*}
for $2\leq k\leq n$ and $a_{k}\in \Gamma (V,M^{(i)}\otimes L)$. Therefore,
evaluating at $p$, since $(a^{*}w)_{p}=(x_{1}^{r-1}dx_{1})_{p}$ (where $%
r\geq 1$ is the index of ramification of $X\rightarrow a(X)$ along $B_{0}$),
we have 
\begin{equation*}
(s\wedge w)_{p}=(a_{1}x_{1}^{r}dx_{1}\wedge ...\wedge dx_{n})_{p}=0.
\end{equation*}
Since $p\in B_{0}$ is general, it follows that $B_{0}\in Bs|\omega
_{X}\otimes M^{(i)}\otimes L|$.
\end{proof}

Our second application is to the case of an $n$-dimensional abelian variety $%
A$. Let 
\begin{equation*}
\left\{ \Theta \right\} \in \mathrm{Pic}\left( A\right)
\end{equation*}
be a positive semidefinite line bundle. Let $Z\subset A$ be any subset. Then
we have the following ``division of linear series'' result.

\begin{theorem}
\label{c10.1}Suppose that there exists positive integers $N,k$ and an ample
divisor 
\begin{equation*}
D=\sum\nolimits_{j\in J_{A}}b_{j}D_{j}\in \Bbb{P}\left( H^{0}\left( A;\Theta
^{N}\right) \right)
\end{equation*}
($D_{j}$ distinct, irreducible) 
such that, for all $j\in J_{A}$, 
\begin{equation*}
nb_{j}<k.
\end{equation*}
Then, for all $L\in Pic^{0}(X)$, there is a section of 
\begin{equation*}
H^{0}\left( \Theta ^{\lceil \frac{nN}{k}\rceil }\otimes \mathcal{I}\left( 
\frac{n}{k}D\right) \otimes L\right)
\end{equation*}
where $\mathcal{I}\left( \frac{n}{k}D\right) $ is the multiplier ideal sheaf
corresponding to the $\Bbb{Q}$ divisor $\frac{n}{k}D$ (see the proof for the
definition). In particular if 
\begin{equation*}
mult_{p}(D)\geq k
\end{equation*}
for all $p\in Z$, then, for all $L\in Pic^{0}(X)$, there is a section of $%
H^{0}\left( \Theta ^{\lceil \frac{nN}{k}\rceil }\otimes L\right) $ which
vanishes on $Z$.
\end{theorem}

\begin{proof}
Following \cite{E}, consider a log resolution 
\begin{equation*}
a:X\longrightarrow A
\end{equation*}
of the pair $(A,D)$, that is, $X$ is obtained from $A$ by a sequence of
blowing-ups with smooth centers lying over $D$ and 
\begin{equation*}
B=\sum\nolimits_{j\in J_{X}}b_{j}B_{j}=a^{*}(D)
\end{equation*}
is a divisor with simple normal crossings. Notice that 
\begin{equation*}
A=\mathrm{Alb}\left( X\right)
\end{equation*}
and 
\begin{equation*}
J_{A}=\left\{ j\in J_{X}:\left. a\right| _{B_{j}}\ generically\
finite\right\}
\end{equation*}
and, for $j\in J_{A}$ we have $D_{j}:=a\left( B_{j}\right) $.

We wish to define the multiplier ideal sheaf ${\mathcal{I}}(\frac{n}{k}D)$
associated to the $\Bbb{Q}$-divisor $\frac{n}{k}D$. Define 
\begin{equation*}
{\mathcal{I}}(\frac{n}{k}D):=a_{*}(K_{X/A}-\lfloor \sum\nolimits_{j\in J}%
\frac{nb_{j}}{k}B_{j}\rfloor )\subseteq \mathcal{O}_{A}.
\end{equation*}
Notice that if $mult _p(D)\geq k$ for all $p\in Z$, then 
\begin{equation*}
{\mathcal{I}}(\frac{n}{k}D)\subseteq {\mathcal{I}}_{Z}.
\end{equation*}
Let 
\begin{equation*}
{\mathcal{F}}:=\Theta ^{\lceil \frac{nN}{k}\rceil }\otimes {\mathcal{I}}(%
\frac{n}{k}D)
\end{equation*}
So we must show that for all $L\in Pic ^0(A)$, 
\begin{equation*}
H^{0}\left( A;{\mathcal{F}} \otimes L\right) \neq 0.
\end{equation*}

Case 1: $\frac{nN}{k}<\lceil \frac{nN}{k}\rceil $.

Then 
\begin{equation*}
\Theta ^{\lceil \frac{nN}{k}\rceil }\otimes \mathcal{O}_{A}(-\frac{n}{k}%
D)\equiv \left( \lceil \frac{nN}{k}\rceil -\frac{nN}{k}\right) \Theta
\end{equation*}
is an ample $\Bbb{Q}$-divisor. Therefore by \cite{E} we have that 
\begin{equation*}
h^{i}\left( {\mathcal{F}}\otimes L\right) =0
\end{equation*}
for all $i>0$ and $\left\{ L\right\} \in {\text{\textrm{Pic}}}^{0}(A)$. The
quantity 
\begin{equation*}
h^{0}({\mathcal{F}}\otimes L)=\chi ({\mathcal{F}}\otimes L)=\chi ({\mathcal{F%
}})\geq 0
\end{equation*}
is constant. It is impossible that $\chi ({\mathcal{F}})=0$, since then $%
h^{i}({\mathcal{F}}\otimes L)=0$ for all $i\geq 0$ and $\left\{ L\right\}
\in \mathrm{Pic}^{0}(A)$ so that, by \cite{M}, ${\mathcal{F}}$ would have to
be the zero sheaf.

%

Case 2: $\frac{nN}{k}=\lceil \frac{nN}{k}\rceil .$

We have $k>nb_j\geq n$. Let 
\begin{equation*}
i:=\frac{nN}{k}<N.
\end{equation*}
We wish to apply the results of \S 8 where 
\begin{equation*}
M=a^{*}\Theta
\end{equation*}
so that 
\begin{eqnarray*}
M^{\left( i\right) } &=&\left( a^{*}\Theta ^{i}\right) \otimes \mathcal{O}
_{X}(-\sum \lfloor \frac{ib_{j}}{N}\rfloor B_{j}) \\
&=&\left( a^{*}\Theta ^{i}\right) \otimes \mathcal{O}_{X}(-\sum \lfloor 
\frac{nb_{j}}{k}\rfloor B_{j})
\end{eqnarray*}
and therefore 
\begin{equation*}
{\mathcal{F}}=a_{*}\left( \omega _{X}\otimes M^{\left( i\right) }\right) .
\end{equation*}
So, by Proposition 1.4 of \cite{E} and the Leray spectral sequence for $%
a_{*} $, 
\begin{equation*}
H^{q}\left( X;\omega _{X}\otimes M^{\left( i\right) }\otimes a^{*}L\right)
=H^{q}\left( A;{\mathcal{F}}\otimes L\right)
\end{equation*}

So it follows from \cite{H} or Theorem \ref{c8.4}ii) that 
\begin{equation*}
V^{i}(\mathcal{F})\supseteq V^{i+1}(\mathcal{F}).
\end{equation*}
If $V^{0}(\mathcal{F})=\emptyset $ then $V^{i}(\mathcal{F})=\emptyset $ for
all $i\geq 0$ and $\chi ({\mathcal{F}})=0$ and we have a contradiction as
above. So, let $S$ be an irreducible component of $V^{0}(\mathcal{F})$, and
denote by 
\begin{equation*}
d:A\longrightarrow S^{\vee }
\end{equation*}
the corresponding map of abelian varieties. We complete the proof by showing
that 
\begin{equation}
S=Pic^{0}\left( A\right) .  \label{spec}
\end{equation}
To prove $\left( \ref{spec}\right) $ suppose that 
\begin{equation*}
0\leq \dim S<n.
\end{equation*}
Then $J_{A}$ is the union of two disjoint subsets 
\begin{equation*}
J_{A}=J^{\prime }\cup J_{S^{\vee }}
\end{equation*}
where $J^{\prime }$ consists of those $j\in J_{A}$ such that 
\begin{equation*}
d\left( D_{j}\right) =S^{\vee }.
\end{equation*}
Notice that 
\begin{equation*}
J^{\prime }\neq \emptyset
\end{equation*}
since 
\begin{equation*}
\sum\nolimits_{j\in J_{S^{\vee }}}b_{j}D_{j}
\end{equation*}
is not ample.

Now by hypothesis 
\begin{equation*}
k>nb_{j}
\end{equation*}
so that it is impossible that 
\begin{equation*}
N\mid ib_{j},
\end{equation*}
that is, that 
\begin{equation*}
k\mid nb_{j}.
\end{equation*}
Thus by Proposition \ref{c9.1} for each $\left\{ L\right\} \in S$ and $j\in
J^{\prime }$ all sections of ${\mathcal{F}}\otimes L$ vanish on $B_{j}$.
Thus we have that 
\begin{equation}
H^{0}\left( \Theta ^{\frac{nN}{k}}\left( -\sum\nolimits_{j\in J^{\prime
}}D_{j}\right) \otimes L\right) =H^{0}\left( \left( a^{*}\Theta ^{\frac{nN}{k%
}}\right) \left( -\sum\nolimits_{j\in J^{\prime }}B_{j}\right) \otimes
a^{*}L\right) \neq 0  \label{9.9'}
\end{equation}
for each $\left\{ L\right\} \in S$. So computing with numerical equivalence
classes on $X$ we have 
\begin{equation}
\begin{array}{l}
\frac{nN}{k}\left\{ a^{*}\Theta \right\} -\sum\nolimits_{j\in J^{\prime
}}\left\{ B_{j}\right\} \\ 
\equiv \sum\nolimits_{j\in J_{X}}\frac{nb_{j}}{k}\left\{ B_{j}\right\}
-\sum\nolimits_{j\in J^{\prime }}\left\{ B_{j}\right\} \\ 
=\sum\nolimits_{j\in J^{\prime }}(\frac{nb_{j}}{k}-1)\left\{ B_{j}\right\}
+\sum\nolimits_{j\in J_{X}-J^{\prime }}\frac{nb_{j}}{k}\left\{ B_{j}\right\}
.
\end{array}
\label{9.10"}
\end{equation}
Let $\left\{ \Theta _{S}\right\} $ be the Chern class of an ample divisor on 
$S^{\vee }$. Then from $\left( \ref{9.10"}\right) $ have on $A$ that 
\begin{equation}
\begin{array}{r}
\left\{ d^{*}\Theta _{S}\right\} ^{\dim S}\cdot \left\{ \Theta \right\}
^{n-\dim S-1}\cdot (\frac{nN}{k}\left\{ \Theta \right\} -\sum\nolimits_{j\in
J^{\prime }}\left\{ D_{j}\right\} ) \\ 
=\left\{ d^{*}\Theta _{S}\right\} ^{\dim S}\cdot \left\{ \Theta \right\}
^{n-\dim S-1}\cdot \sum\nolimits_{j\in J^{\prime }}(\frac{nb_{j}}{k}
-1)\left\{ D_{j}\right\} .
\end{array}
\label{9.11}
\end{equation}
By $\left( \ref{9.9'}\right) $ the expression before the equals sign in $%
\left( \ref{9.11}\right) $ is non-negative. On the other hand the expression
after the equals sign is negative since since $J^{\prime }\neq \emptyset $
and $\frac{nb_{j}}{k}-1<0$ for every $j\in J^{\prime }$.
\end{proof}

A well known result of Esnault and Viehweg (\cite{E} Proposition 5.4 or \cite
{EV2} Theorem 2), concerning zeroes of polynomial equations, states that if
there exists an hypersurface $D$ in $\Bbb{P}^{n}$ of degree $N$ such that $%
mult_{p}(D)\geq k$ for all $p\in Z$, then there is a hypersurface of degree $%
\lfloor \frac{nN}{k}\rfloor $ that contains $Z$. Chudnosky conjectures that
under the above hypothesis, there exists a hypersurface of degree $\lfloor 
\frac{nN}{k}\rfloor -n+1$ that contains $Z$. With the same techniques, one
can obtain similar results for other varieties. Using Theorem \ref{c10.1} we
are able to prove a sharp bound in the case of abelian varieties.

\begin{corollary}
\label{c10.3}Suppose that there exists a divisor 
\begin{equation*}
D=\sum\nolimits_{j\in J_{A}}b_{j}D_{j}\in \Bbb{P}\left( H^{0}\left( A;\Theta
^{N}\right) \right)
\end{equation*}
($D_{j}$ distinct, irreducible) such that 
\begin{equation*}
mult_{p}(D)\geq k
\end{equation*}
for all $p\in Z$. 
Then, there is a section of 
\begin{equation*}
H^{0}\left( \Theta ^{\lceil \frac{nN}{k}\rceil }\right)
\end{equation*}
which vanishes on $Z$.
\end{corollary}

\begin{proof}
If $k\leq n$, the assertion is clear. We may therefore assume that $k>n$. We
may also assume that $n|k$, since if this is not the case, it suffices to
replace $k$ by $nk$, $N$ by $nN$ and $D$ by $nD$.

%
%
Let $N^{\prime }=\frac{k}{n}$. Define 
\begin{equation*}
\Delta :=\sum \lfloor \frac{b_{j}}{N^{\prime }}\rfloor D_{j},
\end{equation*}
\begin{equation*}
\hat{D}:=D-N^{\prime }\Delta =\sum {b}_{j}^{\prime }D_{j}.
\end{equation*}
Therefore, we have that $\frac{n{b^{\prime }}_{j}}{k}<1$. Let $%
Z^{0}:=\left\{ z\in Z:z\in Supp\left( \Delta \right) \right\} $ and $\hat{Z}
=\left( Z-Z^{0}\right) $. It follows for $z\in \hat{Z}$ that 
\begin{equation*}
mult_{z}\hat{D}=mult_{z}D\geq k.
\end{equation*}
Now there is an abelian quotient $A^{\prime }$ of the abelian variety $A$
such that $\Theta \left( -\Delta \right) $ is the pull-back of an ample
divisor 
\begin{equation*}
\Theta ^{\prime }
\end{equation*}
and $\hat{D}$ is the pull-back of a divisor $D^{\prime }\in |N^{\prime
}\Theta ^{\prime }|$ on $A^{\prime }$. Let $Z^{\prime }$ denote the image of 
$\hat{Z}$ in $A^{\prime }.$ Now one can apply the Theorem \ref{c10.1} with 
\begin{equation*}
\Theta =\Theta ^{\prime },N=N^{\prime }=k/n,k=k^{\prime },D=D^{\prime
},Z=Z^{\prime }
\end{equation*}
to obtain a section of $H^{0}(\Theta ^{\prime })=H^{0}(\Theta ^{\frac{nN}{k}
}-\Delta )$ vanishing on $Z^{\prime }$. Since $\Delta $ contains $Z^{0}$,
the theorem follows.
\end{proof}

In \cite{EL} Ein and Lazarsfeld prove that if $(A,\Theta )$ is a principally
polarized abelian variety and $\Theta $ is irreducible, then $(A,\Theta )$
is log terminal (i.e. for all $0\leq \epsilon <1$, the multiplier ideal
sheaf $\mathcal{I}\left( ({1-\epsilon })\Theta \right) $ is trivial).
Moreover if $D\in |N\Theta |$, then the pair $(A,\frac{1}{N}D)$ is log
canonical (i.e. for all $0<\epsilon <1$, the multiplier ideal sheaf $%
\mathcal{I}\left( \frac{1-\epsilon }{N}D\right) $ is trivial). From Theorem 
\ref{c10.1} it is possible to recover the following closely related result
first proved in \cite{H}:

\begin{corollary}
\label{c9.4}Let $(A,\Theta )$ be a principally polarized abelian variety. If 
$N\geq 1$ and $D\in \left| N\Theta \right| $ is such that $\left\lfloor 
\frac{1}{N}D\right\rfloor =0$, then the pair $\left( A,\frac{1}{N}D\right) $
is log terminal.
\end{corollary}

In particular, for any $p\in A$ and a divisor $D$ as above, one has 
\begin{equation*}
mult_{p}D<N\dim A.
\end{equation*}

\begin{proof}
By Theorem \ref{c10.1} (with $k=N\dim A$), for all $L\in Pic^{0}(A)$, one
has $H^{0}(\Theta \otimes \mathcal{I}(\frac{1}{N}D)\otimes L)\ne 0$. It
follows that $\Theta $ vanishes on all the translates of the cosupport of $%
\mathcal{I}(\frac{1}{N}D)$. The only way this can occur is if the cosupport
of $\mathcal{I}(\frac{1}{N}D)$ is empty, i.e. if $\mathcal{I}(\frac{1}{N}D)=%
\mathcal{O}_{A}$.
\end{proof}

\section{Relative statements}

Let $f:X\rightarrow X^{\prime }$ be a surjective map of projective
varieties, $X$ smooth. As before let 
\begin{equation*}
a:X\rightarrow \mathrm{Alb}\left( X\right)
\end{equation*}
and suppose that we have a commutative diagram 
\begin{equation*}
\begin{array}{ccc}
X & \overset{a}{\longrightarrow } & \mathrm{Alb}\left( X\right) \\ 
\downarrow ^{f} &  & \downarrow ^{g} \\ 
X^{\prime } & \overset{a^{\prime }}{\longrightarrow } & A^{\prime }
\end{array}
\end{equation*}
where $g$ is a morphism of abelian varieties. Let 
\begin{equation*}
P^{\prime }\rightarrow A^{\prime }\times \mathrm{Pic}^{0}\left( A^{\prime
}\right)
\end{equation*}
be the Poincar\'{e} bundle. Let 
\begin{equation*}
\Delta
\end{equation*}
denote a formal neighborhood in $g^{\vee }\left( \mathrm{Pic}^{0}\left(
A^{\prime }\right) \right) \subseteq \mathrm{Pic}^{0}\left( X\right) $ of a
point $\left\{ L_{0}\right\} \in g^{\vee }\left( \mathrm{Pic}^{0}\left(
A^{\prime }\right) \right) $ and let $\zeta ^{\prime }$ denote the linear
local coordinate for $g^{\vee }\left( \mathrm{Pic}^{0}\left( A^{\prime
}\right) \right) $ with $\left\{ L_{0}\right\} $ given by $\zeta ^{\prime
}=0 $.

In \cite{Ko2} Koll\'{a}r shows that 
\begin{equation}
R^{\cdot }f_{*}\omega _{X}\cong \sum\nolimits_{i}R^{i}f_{*}\omega _{X}[-i].
\label{10.1}
\end{equation}
In particular 
\begin{equation*}
h^{p}(X,\omega _{X})=\sum\nolimits_{i}h^{i}(X^{\prime },R^{p-i}f_{*}\omega
_{X}).
\end{equation*}
So by the projection formula 
\begin{eqnarray*}
&&R^{\cdot }(\pi \circ f\times id_{\Delta })_{*}\left( \omega _{X}\otimes
\left( g\circ a\right) ^{*}P^{\prime }\right) \\
&=&R^{\cdot }\pi _{*}R^{\cdot }(f\times id_{\Delta })_{*}\left( \omega
_{X}\otimes \left( g\circ a\right) ^{*}P^{\prime }\right) \\
&=&R^{\cdot }\pi _{*}\left( R^{\cdot }f_{*}(\omega _{X})\otimes \left( {a}
^{\prime }\right) ^{*}P^{\prime }\right) \\
&=&\sum R^{\cdot }\pi _{*}\left( R^{i}f_{*}(\omega _{X})\otimes \left( {a}
^{\prime }\right) ^{*}P^{\prime }\right) [-i].
\end{eqnarray*}

Now 
\begin{equation*}
\left. \beta \right| _{\Delta }=\beta ^{\prime }\circ g
\end{equation*}
for some linear functional $\beta ^{\prime }\in \mathrm{Hom}_{\Bbb{R}}\left(
H_{1}\left( A^{\prime };\Bbb{Z}\right) ,\Bbb{R}\right) $. Referring to
Theorem \ref{c7.1} it follows that 
\begin{eqnarray*}
&&\sum R^{q-i}\pi _{*}\left( R^{i}f_{*}\omega _{X}\otimes L_{0}\otimes
\left( {a}^{\prime }\right) ^{*}P^{\prime }|_{\Delta _{\frak{A}}}\right) \\
&=&R^{q}(\pi \circ f\times id_{\Delta })_{*}\left( \omega _{X}\otimes
g^{*}L_{0}\otimes a^{*}\left( g^{*}P^{\prime }\right) |_{\Delta _{\frak{A}
}}\right) \\
&\cong &\mathcal{H}^{q}\left( H^{*}(\omega _{X}\otimes g^{*}L_{0})\otimes 
\frac{\Bbb{C}\left[ \left[ \zeta ^{\prime }\right] \right] }{\frak{A}},%
\overline{\partial }\beta \wedge \right) \\
&=&\sum \mathcal{H}^{q-i}\left( H^{*}(R^{i}f_{*}\omega _{X}\otimes
L_{0})\otimes \frac{\Bbb{C}\left[ \left[ \zeta ^{\prime }\right] \right] }{%
\frak{A}},\overline{\partial }\beta ^{\prime }\wedge \right) .
\end{eqnarray*}
But, referring to the grading in \S 6, the associated isomorphism 
\begin{equation*}
\begin{array}{r}
Gr\left( R^{q}(\pi \circ f\times id_{\Delta })_{*}\left( \omega _{X}\otimes
g^{*}L_{0}\otimes a^{*}\left( g^{*}P^{\prime }\right) |_{\Delta _{\frak{A}
}}\right) \right) \\ 
\cong \mathcal{H}^{q}\left( H^{*}(\omega _{X}\otimes g^{*}L_{0})\otimes 
\frac{\Bbb{C}\left[ \left[ \zeta ^{\prime }\right] \right] }{\frak{A}},%
\overline{\partial }\beta \wedge \right)
\end{array}
\end{equation*}
is natural and sends each summand 
\begin{equation*}
Gr\left( R^{q-i}\pi _{*}\left( R^{i}f_{*}\omega _{X}\otimes L_{0}\otimes
\left( {a}^{\prime }\right) ^{*}P^{\prime }|_{\Delta _{\frak{A}}}\right)
\right)
\end{equation*}
to the corresponding summand 
\begin{equation*}
\mathcal{H}^{q-i}\left( H^{*}(R^{i}f_{*}\omega _{X}\otimes L_{0})\otimes 
\frac{\Bbb{C}\left[ \left[ \zeta ^{\prime }\right] \right] }{\frak{A}},%
\overline{\partial }\beta ^{\prime }\wedge \right) .
\end{equation*}
Thus:

\begin{theorem}
\begin{equation*}
Gr\left( R^{j}\pi _{*}\left( R^{i}f_{*}\omega _{X}\otimes L_{0}\otimes
\left( {a}^{\prime }\right) ^{*}P^{\prime }|_{\Delta _{\frak{A}}}\right)
\right) \cong \mathcal{H}^{j}\left( H^{*}(R^{i}f_{*}\omega _{X}\otimes
L_{0})\otimes \frac{\Bbb{C}\left[ \left[ \zeta ^{\prime }\right] \right] }{%
\frak{A}},\overline{\partial }\beta ^{\prime }\wedge \right) .
\end{equation*}
\end{theorem}

\begin{proof}
Kollar's celebrated result $\left( \ref{10.1}\right) $ says that the Leray
spectral sequence for $f_{*}$ degenerates at $E_{2}$. But the filtration in
the Leray spectral sequence is compatible with that of the spectral sequence
defined in the proof of Theorem \ref{c6.1}i).
\end{proof}

\section{Appendix: Existence of compatible trivializations of line bundles}

\begin{proposition}
Given a transversely holomorphic trivialization 
\begin{equation*}
F_{\sigma }:X\overset{\left( \sigma ,\pi \right) }{\longrightarrow }%
X_{0}\times \Delta
\end{equation*}
of a deformation $X/\Delta $ and given a holomorphic line bundle $L/X$, a
compatible trivialization of line bundles $\left( \ref{2.3.1}\right) $
always exists.
\end{proposition}

\begin{proof}
Let $\left\{ W\right\} $ be a covering of $X$ by coordinate disks and $%
\left\{ W_{0}\right\} $ the restriction of this covering to $X_{0}$. We
construct a $C^{\infty }$ partition-of-unity $\left\{ \rho _{W_{0}}\right\} $
subordinate to the induced covering of $X_{0}$. Recall that $L$ is given
with respect to the trivialization $\sigma $ by holomorphic local patching
data 
\begin{eqnarray*}
g^{_{WW^{\prime }}}\left( x\right) &=&\sum g_{i}^{_{WW^{\prime }}}\left(
x_{0}\right) t^{i} \\
&=&g^{W_{0}W_{0}^{\prime }}\left( x_{0}\right) \exp \left(
\sum\nolimits_{j>0}a_{j}^{_{WW^{\prime }}}\left( x_{0}\right) t^{j}\right)
\end{eqnarray*}
where $x_{0}=\sigma \left( x\right) $ and 
\begin{equation*}
\sum\nolimits_{j>0}a_{j}^{_{WW^{\prime }}}\left( x_{0}\right) t^{j}=\log 
\frac{g^{_{WW^{\prime }}}\left( x\right) }{g^{W_{0}W_{0}^{\prime }}\left(
x_{0}\right) }.
\end{equation*}
Notice that, if $V,W,$ and $W^{\prime }$ are three open sets of the cover
which have non-empty intersection, then, for all $j>0$, 
\begin{equation*}
a_{j}^{_{VW}}+a_{j}^{_{WW^{\prime }}}=a_{j}^{_{VW^{\prime }}}.
\end{equation*}
Define the mapping 
\begin{equation*}
L\rightarrow L_{0}
\end{equation*}
over $W_{0}\times \Delta $ by 
\begin{equation}
\left( x,v\right) \mapsto \left( x_{0},\exp \left( \sum\nolimits_{W^{\prime
}}\rho _{W_{0}^{\prime }}\left( x_{0}\right) \left(
\sum\nolimits_{j>0}a_{j}^{_{WW^{\prime }}}\left( x_{0}\right) t^{j}\right)
\right) \cdot v\right) .  \label{2.3.0}
\end{equation}
This map is well defined since, over $V\cap W$ we have 
\begin{equation*}
g^{VW}\left( x\right) =g^{V_{0}W_{0}}\left( x_{0}\right) \exp \left(
\sum\nolimits_{j>0}a_{j}^{_{VW}}\left( x_{0}\right) t^{j}\right)
\end{equation*}
and so 
\begin{eqnarray*}
&&g^{_{VW}}\left( x\right) \cdot \exp \left( \sum\nolimits_{W^{\prime }}\rho
_{W_{0}^{\prime }}\left( x_{0}\right) \left(
\sum\nolimits_{j>0}a_{j}^{_{WW^{\prime }}}\left( x_{0}\right) t^{j}\right)
\right) \\
&=&g^{V_{0}W_{0}}\left( x_{0}\right) \exp \left(
\sum\nolimits_{j>0}a_{j}^{_{VW}}\left( x_{0}\right) t^{j}\right) \cdot \exp
\left( \sum\nolimits_{W^{\prime }}\rho _{W_{0}^{\prime }}\left( x_{0}\right)
\left( \sum\nolimits_{j>0}a_{j}^{_{WW^{\prime }}}\left( x_{0}\right)
t^{j}\right) \right) \\
&=&g^{V_{0}W_{0}}\left( x_{0}\right) \exp \left( \sum\nolimits_{W^{\prime
}}\rho _{W_{0}^{\prime }}\left( x_{0}\right) \sum\nolimits_{j>0}\left(
a_{j}^{_{VW}}+a_{j}^{_{WW^{\prime }}}\right) \left( x_{0}\right) t^{j}\right)
\\
&=&g^{V_{0}W_{0}}\left( x_{0}\right) \exp \left( \sum\nolimits_{W^{\prime
}}\rho _{W_{0}^{\prime }}\left( x_{0}\right)
\sum\nolimits_{j>0}a_{j}^{_{VW^{\prime }}}\left( x_{0}\right) t^{j}\right) .
\end{eqnarray*}
\end{proof}

\end{document}